\documentclass[preprint]{elsarticle}

\usepackage{amsmath}
\usepackage{amssymb}
\usepackage{amstext}
\usepackage{color}

\definecolor{purple}{rgb}{0.65, 0, 1}
\definecolor{orange}{rgb}{1,.5,0}
\definecolor{juice}{rgb}{.9,.43,.26}

\usepackage{graphicx}
\usepackage{hhline}
\usepackage{psfrag}
\usepackage{float}
\usepackage[dvips]{epsfig}
\usepackage{subfig}

\newcommand{\R}{\mathrm{I\hspace{-0.5ex}R}}

\newcommand{\dd}{\mathrm{d}}

\newtheorem{theorem}{Theorem}
\newtheorem{remark}{Remark}

\newdefinition{algorithm}{Algorithm}
\newdefinition{algorithm2}{Algorithm}

\newtheorem{corollary}{Corollary}
\newtheorem{proposition}{Proposition}
\newtheorem{definition}{Definition}
\newtheorem{lemma}{Lemma}
\newproof{proof}{Proof}
\newcommand{\rr}{\mathbb{R}}
\newcommand{\rrn}{\mathbb{R}^{n}}

\renewcommand{\div}{\operatorname{div}}

\begin{document}
\begin{frontmatter}

\title{A new interaction potential for swarming models
}
\author[sm]{J. A. Carrillo 
}
\ead{carrillo@imperial.ac.uk}
\author[sm]{S. Martin \corref{cor1}
}
\ead{stephan.martin@imperial.ac.uk}
\author[vp]{V. Panferov
}
\ead{vladislav.panferov@csun.edu}

\address[sm]{Department of Mathematics, Imperial College London, London SW7 2AZ, UK \\
Phone: +44- (0) 20-759-48396, Fax: +44- (0) 20-759-48517}
\address[vp]{Department of Mathematics, California State University Northridge \\
Northridge, CA 91330-8313,USA }

\cortext[cor1]{Corresponding author}

\begin{abstract}
{ We consider a self-propelled particle system which has been used
to describe certain types of collective motion of animals, such as
fish schools and bird flocks. Interactions between particles are
specified by means of a pairwise potential, repulsive at short
ranges and attractive at longer ranges. The exponentially decaying
Morse potential is a typical choice, and is known to reproduce
certain types of collective motion observed in nature,
particularly aligned flocks and rotating mills. We introduce a
class of interaction potentials, that we call Quasi-Morse, for
which flock and rotating mills states are also observed
numerically, however in that case the corresponding macroscopic
equations allow for explicit solutions in terms of special
functions, with coefficients that can be obtained numerically
without solving the particle evolution. We compare the obtained
solutions with long-time dynamics of the particle systems and find
a close agreement for several types of flock and mill solutions. }
\end{abstract}


\end{frontmatter}


\section{Introduction}
Emerging behaviors in interacting particle systems
have received a lot of attention in research in recent years.
Topics range from diverse fields of applications such as animal
collective behavior, traffic, crowd dynamics and crystallization.
Self-organization in the absence of leaders has been reported in
several species which coordinate their movement (swarming), and
several models have been proposed for their explanation
\cite{parrish,mogilner,BDT,camazine,couzin,CKJRF}.

Many of these models are based on zones in which some of 3 basic
effects are included: short-range repulsion, long-range
attraction, and alignment. These 3-zone basic
descriptions have been very popular for modeling fish schools
\cite{HW,KH,HK,BTTYB,BEBSVPSS}, starlings \cite{HCH}, or ducks
\cite{LLE,LLE2}. The main modelling issues are if some or all of
these effects between agents have to be included and if so, how to
incorporate them. Many basic swarming models rely on averaged
spatial distance or orientation interactions while recent
biological studies point out the importance of nearest-neighbor
interactions \cite{rome} or anisotropic communication
\cite{KTIHC}. Mathematicians have started in recent years to
attack one of the most striking features of these {\it simple
looking} models: the diversity of swarming states, also called
patterns in the biology community, their emergence and stability.

The individual level description of these phenomena
leads to certain particle systems, called Individual Based Models
(IBMs), with some common aspects. Typically, the attraction-repulsion is modeled
through pairwise effective potentials depending on the
distance between individuals. An asymptotic speed for particles is
imposed either by working in the constrained set of a sphere in
velocity space \cite{vicek,chate,DM1} or by adding a term of
balance between self-propulsion and friction which effectively
fixes the speed to a limiting value for large times
\cite{LR,DCBC}. In this work, we will not include any alignment
mechanism. We refer to \cite{review} for a survey on results
related to kinetic modeling in swarming.

In Section 2 we will review some of these IBMs, and discuss the
appearance of two main swarming patterns: mills and flocks. These
patterns are easily observed in particle simulations
\cite{DCBC,review2} and reported in detail for certain particular
potentials, the so-called Morse potentials. We will give a precise
definition of flocks and mills as solutions of the kinetic
equation associated to the particle systems.
Finding the spatial shape of flocks and mills has
been numerically reported in the literature but obtaining
analytical results on them has only been done in one dimension for
the Morse potential in \cite{BT2}.

In this work, we generalize the strategy in
\cite{BT2} proposing a new interaction potential, that we call
Quasi-Morse, to replace the Morse potential. The Quasi-Morse
potential coincides with the Morse potential in one dimension and
we will show that it is a suitable extension of the Morse
potential in $n=2,3$. Section 3 introduces Quasi-Morse potentials
as fundamental solutions of certain linear PDEs. We will first
show that the Quasi-Morse potentials are biologically relevant in
essentially the same parameter range as the Morse potentials.
Second, we make use of their particular structure to show in our
main theorem that flock and mill solutions can be expressed as
almost explicit linear combinations of special functions.

Finally, Section 4 is devoted to propose an algorithm to compute
the scalar coefficients in the expansion of the flock and mill
patterns in terms of the basis functions associated with the
Quasi-Morse PDE operators. The strategy uses ideas of constrained
optimization methods. We finally compare the results for flocks in
2D and 3D and mills in 2D to particle simulations showing a good
agreement. 
As a conclusion, we demonstrate that the proposed
Quasi-Morse potentials are a very good alternative to Morse
potentials as they share many of their features in the natural
parameter range, and at the same time enable explicit computation
of the macroscopic density profiles up to numerically determined
constants.

\section{Swarming: Models \& Patterns}
\label{section2}

We will consider a simple second order model for swarming analyzed
in \cite{DCBC} consisting of the attraction-repulsion of $N$
interacting self-propelled particles located at
$x_{i}\in\R^{n}$ with velocities $v_{i}\in\R^{n}$ in a host medium
with friction, with $n=1,2,3$. Friction is modeled by Rayleigh's
law and as a result, an asymptotic speed for the individuals is
fixed by the compensation of friction and self-propulsion. More
precisely, the time evolution is governed by the
equations of motion
\begin{flalign}\label{micro}
\begin{split}
   \frac{dx_i}{dt} &= v_i \,, \\
   \frac{dv_i}{dt} &= \alpha v_i  - \beta v_i \vert v_i\vert^2
   - \nabla_{x_i} \sum_{i \neq j} W(x_i-x_j)\,,
\end{split}
\end{flalign}
where $W$ is a pairwise interaction potential and $\alpha, \beta$ are
effective values for propulsion and friction forces, see
\cite{LR,DCBC,CDMBC,CHDB} for more discussion. The interaction
potential $W: \rr^n \times \rr^n \rightarrow \rr $ is assumed to
be radially symmetric: $W(x) = U(|x|)$, $x\in\rr^n$. The typical
asymptotic speed of the individuals is $\sqrt{\alpha/\beta}$.
The Morse potential is defined by taking
\begin{align*}
    U(r)  = -C_A e^{-r/l_A}  + C_R e^{-r/l_R} ,
\end{align*}
where $C_A$, $C_R$ are the attractive and repulsive strengths, and
$l_A$, $l_R$ are their respective length scales. We set
$V(r)=-\exp(-r/l_A)$, $C=C_R/C_A$, and $l=l_R/l_A$ to obtain
\[
 U(r) = C_A \left[ V(r) - C V\left(\frac{r}{l}\right) \right].
\]

The choice of this potential is motivated in \cite{DCBC} for being
one of the simplest choices of integrable potentials with easily
computable conditions to distinguish the relevant parameters in
biological swarms. In fact, it is straightforward to check that in
the range $C>1$ and $l<1$ the potential $U(r)$ is short-range
repulsive and long-range attractive with a unique minimum defining
a typical distance between particles. Moreover, in this regime the sign of the
integral of the potential:
\begin{equation}\label{morse}
\mathcal{U}:=\int_0^\infty W(x)\,dx = \mathcal{V} (1-Cl^n) \qquad
\mbox{with } \mathcal{V}:=\int_0^\infty V(r) r^{n-1}\,dr<0\,,
\end{equation}
gives a criterion to distinguish between the so-called
 H-stable and catastrophic regimes. This condition reads as $C l^n - 1 <0$ for the catastrophic case in any dimension $n$, see \cite{DCBC,rue}. This property of the
potential is important since it is related to the typical patterns emerging in such systems, as classified in
\cite{DCBC}.

Flocks, where particles tend to form groups, moving with
the same velocity, and milling solutions, where rotatory states are formed
are of particular interest and are observed in particle and hydrodynamic
simulations \cite{DCBC,CKMT} in $n=2$. Actually, they typically emerge in the large time behavior of the system of
particles \eqref{micro} in the catastrophic regime $C l^2 <1$ with
$C>1$ and $l<1$. In the same range of parameters, randomly chosen
initial data lead also to other patterns such as double mills and
flocks \cite{DCBC,CDP}. However mills are not observed in the
H-stable regime $C l^2>1$ with $C>1$ and $l<1$ while flocks do.

Assuming the weak coupling scaling~\cite{dobru,neunzert,BH,spohn2}
in which the range of interaction is kept fixed and the strength of interaction
is divided proportionally between particles, we pass to the rescaled formulation:
\begin{flalign*}
\begin{split}
  \frac{dx_i}{dt} &= v_i \, ,\\
  \frac{dv_i}{dt} &= v_i  ( \alpha  - \beta \vert v_i\vert^2 )
  - \frac{1}{N} \nabla_{x_i} \sum_{i \neq j} U(\vert x_i - x_j \vert
  )\, .
\end{split}
\end{flalign*}
This system has a well-defined limit as \(N\to\infty\) which can be expressed as
a solution of the corresponding mean-field equation:
\begin{equation}
\partial_{t}f+v\cdot \nabla_{x}f + F[\rho]\cdot \nabla_{v}f +\div\left(\left(\alpha - \beta |v|^{2}\right)v f\right)
=0\,, \label{VP}
\end{equation}
with
\begin{equation*}
\rho(t,x):=\int f(t,x,v)\dd v \,.
\end{equation*}
Here, $f(t,x,v): \rr \times \rrn \times \rrn \rightarrow \rr$ is
the phase-space density, and $\rho(t,x)$ is the averaged
(macroscopic) density. The mean-field interaction is given by
$F[\rho] =-\nabla_{x} W \star \rho$.

The limit \(N\to\infty\) has been established rigorously
for smooth potentials $W\in C^2_b$ in \cite{dobru,neunzert,BH,spohn2}, in
\cite{CCR,BCC} for more general models with and without noise, and
for more general potentials, with possibly singular behavior at
zero, including the Morse potential \eqref{morse} in the recent
result \cite{hj}.

\subsection{Flock and Mill States}

We are interested in computing certain relevant particular
solutions of the Vlasov-like equation for swarming in \eqref{VP}.
In fact, we can formally find mono-kinetic solutions of \eqref{VP}
by inserting the ansatz:
\begin{equation}\label{monokinetic}
f(t,x,v ) = \rho (t,x) \, \delta (v - u(t,x)),
\end{equation}
in the weak formulation of \eqref{VP}. The result in
\cite{CDMBC,CDP} is that $\rho$ and $u$ should satisfy the
following set of hydrodynamic equations:
\begin{equation}
\left\lbrace
\begin{array}{l}
\displaystyle \frac{\partial \rho}{\partial t} +
\div_{x}(\rho u) = 0, \vspace{.3cm}\\
\rho\,\displaystyle \frac{\partial u}{\partial t} + \rho\,(u
\cdot\nabla_x) u = \rho\, (\alpha - \beta |u|^2) u -
\rho\,(\nabla_x W \star \rho).
\end{array}
\right. \label{momentsimpl}
\end{equation}

\begin{definition}
A \emph{flock} is a solution $f_{F}$ of \eqref{VP} of the form
\eqref{monokinetic} with $\rho (t,x)=\rho_{F}(x-tu_{0})$ and
$u(t,x)=u_0$ with $u_0\in \rr^n$ such that
$|u_{0}|=\sqrt{\frac{\alpha}{\beta}}$ and $\rho_F$ a probability
measure in $\rr^n$.
\end{definition}

\

Obviously, flock solutions are determined by their density profile
$\rho_{F}$ and have the structure of traveling waves in the direction
of the velocity vector \(u_0\). It is straightforward to see that the
density of a flock is characterized by the following equation:
\begin{proposition}
The function $f_F\ge 0$ is a flock solution if and only if the macroscopic density \(\rho_F\)
satisfies
\begin{equation}\label{flockchar}
\nabla_{x} W \star \rho_F = 0 \quad \mbox{ on the support of }
\rho_F\, .
\end{equation}
\end{proposition}
There are singular solutions to \eqref{flockchar} obtained by
concentrating all the mass uniformly in a suitable sphere, the
so-called Delta rings \cite{CDP}, whose stability for first order
models has recently been studied in \cite{BCLR} for
certain potentials. Also, there are solutions to \eqref{flockchar}
given by smooth compactly supported densities for combination of
suitable powers in 1D
\cite{FellnerRaoul1,FellnerRaoul2}, for the Morse potential in 1D
\cite{BT2}, and for combination of powers when one of them is the
repulsive Newtonian potential \cite{FHK} in 2D. In fact, the set
of solutions to \eqref{flockchar} can be very complicated even in
one dimension \cite{FellnerRaoul1,FellnerRaoul2,HLF} depending on
the regularity of the potential.

Let us remark that since we assume the radial symmetry of the
potential, one expects that the density of the flocking solutions
to \eqref{flockchar} is radially symmetric as well and that it is
supported in a ball $B(0,R_F)$ with $R_F>0$. This is reinforced by
the fact that the convolution of radial functions is radial, see
Subsection 2.2 for more precise statements.
We will reduce ourselves to find flocking
solutions with radial symmetry in the rest of this paper, that is,
finding $R_F>0$ and a radial density $\rho_F(|x|)$ compactly
supported in $B(0,R_F)$ satisfying
\begin{equation}\label{flockcharrad}
W \star \rho_F = C \quad \mbox{ in } B(0,R_F)\,,
\end{equation}
for some constant $C\in\rr$.

Another interesting type of solutions that spontaneously show up
in particle simulations are {\em  mills}, they correspond
to motion with the velocity field of a point vortex:
\begin{equation}\label{vortex}
u_M(x) = \pm \sqrt{\frac {\alpha} {\beta}} \, \frac
{x^{\perp}}{|x|}\,,
\end{equation}
where $x=(x_1,x_2)$, $x^{\perp}=(-x_2,x_1)$, such that $\rho_M(|x|)$ is a radially
symmetric stationary solution to
\eqref{momentsimpl}.
\begin{definition}
A \emph{mill} is a solution $f_{M}$ of \eqref{VP} of the form:
\begin{equation*}
 f_{M}(t,x,v) = \rho_{M}(x) \, \delta (v-u_M(x))\,,\,
\end{equation*}
with $u_M$ given by \eqref{vortex} and $\rho_M$ radially symmetric.
\end{definition}

As shown in \cite{LR,CDP,CKMT}, mill solutions can
also be characterized as:
\begin{proposition}
$\rho_M(x) $ is a mill density if and only if
\begin{equation*}
\nabla_x \left[ W \star \rho - \frac {\alpha}{\beta} \log
|x|\right]=0, \qquad \text{on the support of } \rho \, .
\end{equation*}
\end{proposition}

As discussed above, one can obtain singular mill solutions by
concentrating all particles in a ring \cite{CDP}. However, we will
search for radial solutions supported in an annulus $B(R_m,R_M)$
with $0<R_m<R_M$, and therefore, mill radial solutions supported
in $B(R_m,R_M)$ are characterized by
\begin{equation}\label{millcharrad}
W \star \rho_M = D + \frac {\alpha}{\beta} \log |x| \quad \mbox{
in } B(R_m,R_M)\, ,
\end{equation}
for some constant $D\in\rr$.  In the following,
the subindex $x$ in differential operators is dropped since we
only deal with $x$-dependent functions.

\subsection{Convolution of radial functions}
\label{sectionradialconvolutions}

Since we want to find particular radial solutions
to flocks \eqref{flockcharrad} and mills \eqref{millcharrad}, we
need suitable expressions of the convolution of two radial
functions in $n=2,3$. Given any radial density $\rho(|x|)$, then
the convolution term rewrites:
$$
 (W\star \rho)(x) = \int_{\rrn} W(x-y)\rho(|y|)\dd y = \int^{\infty}_{0}\int _{\partial
 B(0,1)}W(x-s\omega)\rho(s)s^{n-1} \dd \omega\dd s
$$
which is \emph{not} a convolution in $r=|x|$ anymore, but rather
is given by an integral operator of the following form:
\begin{equation*}
 (W\star \rho) (r) =  \int_{\rr^{+}} \Psi(r,s) \rho(s)\dd s
\end{equation*}
with
\begin{equation*}
\Psi(r,s) = s^{n-1}\int_{\partial B(0,1)} U(|re_1-s\omega| )\dd
\omega\,.
\end{equation*}
Expressing it in polar ($n=2$) or spherical ($n=3$) coordinates,
we get the functions
\begin{equation}\label{conv2d}
 \Psi(r,s)=s \int_{0}^{2\pi} U\left(\sqrt{r^{2}-2rs\cos\theta+s^{2}} \right)
 \dd\theta
\end{equation}
for $n=2$ and
\begin{align}
\Psi(r,s)& =s^2\int^{2\pi}_{0}\!\int^{\pi}_{0}U(|r
e_1-s\omega(\theta,\nu)|)\,\sin\nu\,\dd\nu\,\dd\theta \nonumber\\
&= 2\pi s^2 \int^{\pi}_{0}U\left(\sqrt{r^{2}-2rs\cos\nu+s^{2}}
\right)\,\sin\nu\,\dd\nu\,,\label{conv3d}
\end{align}
with $\omega(\theta,\nu)=(\cos \nu, \sin\nu \cos\theta,\sin\nu
\sin\theta)$ for $n=3$.


\section{Quasi-Morse potentials and their explicit solvability}
\label{section3} In this section, we define Quasi-Morse potentials
for $n=1,2,3$ and discuss their properties. These
Quasi-Morse potentials will yield biologically relevant shapes
similar to the Morse potentials. We show that flock and mill
solutions in the natural parameter range, see Figures
\ref{fig-2dflockparameters} and \ref{fig-3d}(d) for precise
statements, can be computed explicitly up to numerically
determined constants.

\subsection{Definition and comparison}

\begin{definition}\label{defv}
Let $V:\rr^+\rightarrow\rr$ denote the radially symmetric solution
of the $n$-dimensional \emph{screened Poisson equation} $\Delta u
- k^{2} u = \delta_{0}$, for a given $k>0$, that vanishes at
infinity. Let $C,l,\lambda \in\rr$ be further positive parameters.
Then we say that $U(|x|)$ is the \emph{\(n\)-dimensional
Quasi-Morse potential} if
\begin{equation*}
U(r):=\lambda\left(V(r)-C\, V\left(\frac{r}{l}\right)\right)\, .
\end{equation*}
\end{definition}

\begin{figure}
\centering \subfloat[Quasi-Morse
potential]{\includegraphics[keepaspectratio=true,
width=.4\textwidth]{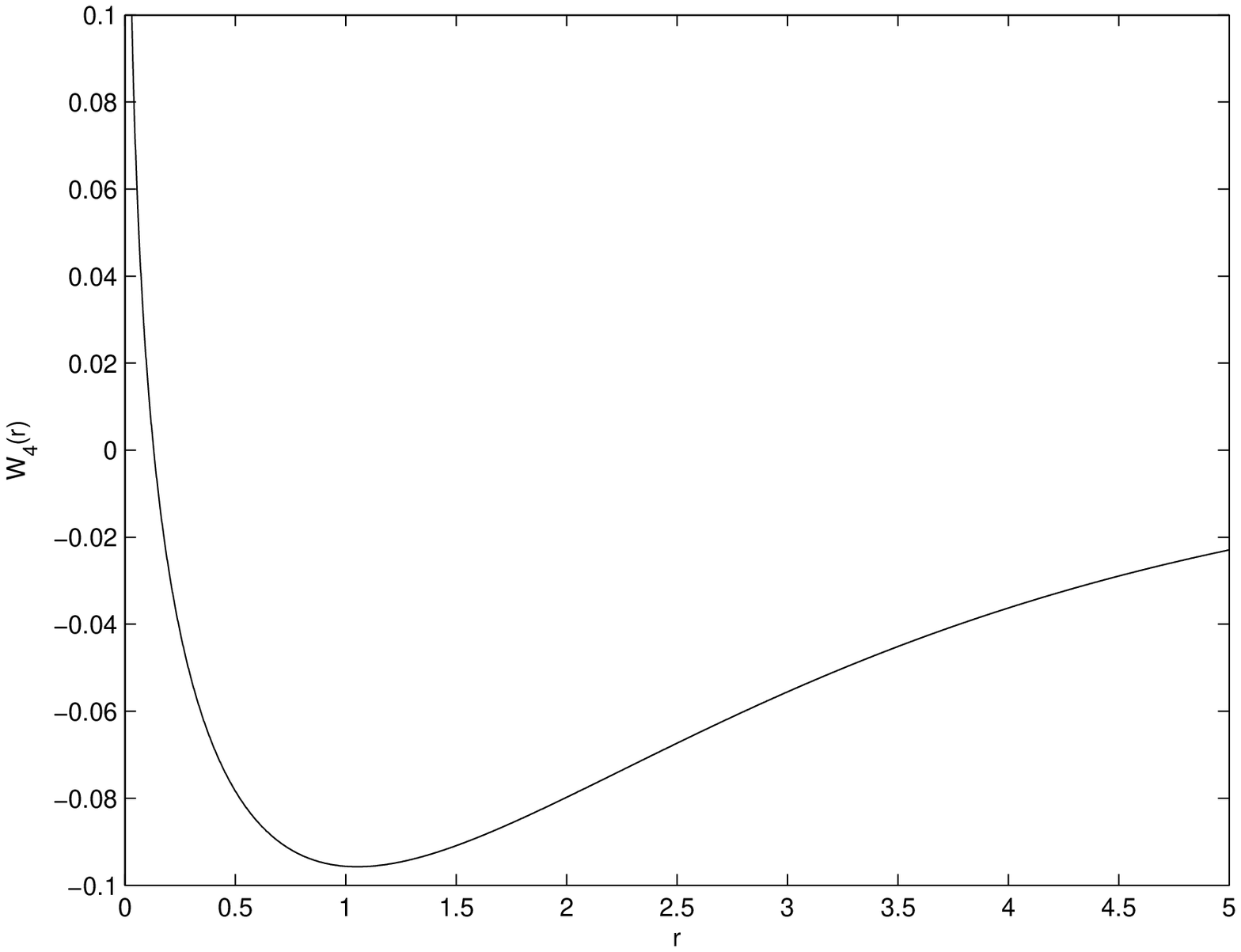}}
\subfloat[Morse potential]{\includegraphics[keepaspectratio=true,
width=.4\textwidth]{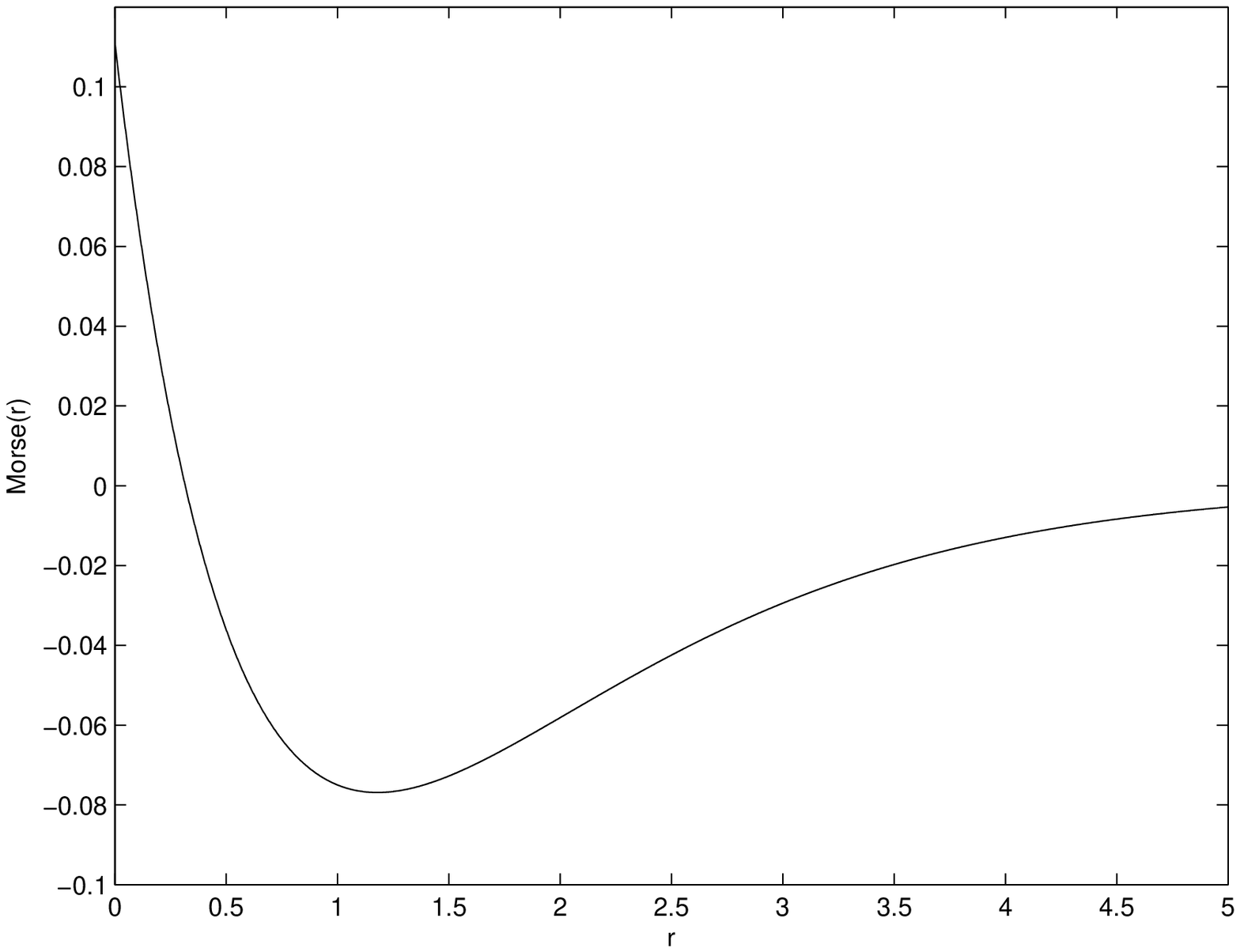}}
\caption{Comparison of potentials: Both yield the biologically
relevant shape of short-range repulsion and long-range attraction
(Quasi-Morse: $n=2, C=\frac{10}{9}, l=0.75, k=\frac{1}{2},
\lambda=4$, Morse: $C=\frac{10}{9}, l=0.75, \lambda=2$).}
\label{fig-comparepotentials}
\end{figure}

Using the radially symmetric ansatz, the screened Poisson equation
reduces to a second-order ordinary differential equation dependent
on the space dimension. For relevant $n=1,2,3$ this ODE possesses two
linearly independent solutions. We therefore have
\begin{corollary}
\label{cor-quasimoorsedef}
Quasi-Morse potentials for $n=1,2,3$ are well-defined and
constructed from the following fundamental solution:
\begin{equation}\label{potdef}
\begin{cases}
n=1: & V(r) =-\frac{1}{k}e^{-kr} \\
n=2: & V(r) = -\frac{1}{2\pi}K_{0}(kr) \\
n=3: & V(r) = -\frac{1}{4\pi}\frac{e^{-kr}}{r}
\end{cases}
\end{equation}
where $K_{0}$ is the modified Bessel function of second kind. For
$n=1$, the Quasi-Morse potential equals the Morse potential.
\end{corollary}
We illustrate the Quasi-Morse potential in comparison to the Morse
potential for $n=2$ with parameters $C=10/9$ and $l=0.75$ in
Figure \ref{fig-comparepotentials}.
Both potentials could be used to model the biologically motivated
interplay between short-range repulsion and long-range attraction,
and there is no clear reason to prefer one over the other. A
significant difference is the behavior at zero, where Morse is
finite and Quasi-Morse is singular though locally integrable for
$n>1$, which are the dimensions we aim to study.
The parameter dependence of catastrophic regimes
is inherited from the Morse potential as summarized in the next
result, whose proof is given in an appendix.

\begin{corollary}\label{cor-biorel}
The function \(U(r)\) has a unique minimum if and only if \(l<1\),
\(Cl^{n-2}>1\). Furthermore, the Quasi-Morse potential \(U(|x|)\)
is catastrophic if $Cl^{n}<1$.
\end{corollary}

\begin{remark}
We first emphasize that the global minimum of $U$
corresponds to the biologically relevant scenario of short-range
repulsion and long-range attraction, as for the standard Morse
potential. Concerning the H-stability of the Quasi-Morse
potentials, we remark that the inverse Fourier transform of $U(r)$
for $k=1$ reads
\[
\check{U}(|\xi|) =
\frac{Cl^{n}-1+l^2(Cl^{n-2}-1)|\xi|^2}{(1+|\xi|^2)(1+l^2\,|\xi|^2)}.
\]
which is positive if  $Cl^{n}>1$ and $Cl^{n-2}>1$. This indicates the 
H-stability, but the criteria developed in \cite{rue} do not apply
directly, since $\check{U}(|\xi|)$ is not integrable in dimensions
$n=2,3$. However, our numerical findings
presented in the following sections will suggest H-stability for
the configurations $l<1, Cl^{n-2}>1, Cl^{n}>1$. This corresponds
to potentials having a unique minimum and a positive
$n-$dimensional integral.
\end{remark}

Next, we mention the influence of the free scaling parameter $k$
and show that any potential shape can be normalized to $k=1$,
see appendix. The following results are given
without proof which follows easily by a change of variables from
the convolution form in radial coordinates \eqref{conv2d} and
\eqref{conv3d} in subsection 2.2.

\begin{corollary}
Let $\rho$ be the flock (resp.\@ mill) solution setting $k=1$ with
support $B(0,R)$ (resp.\@ $B(R_{m},R_{M})$), then the transformed
solution $\tilde{\rho}$ for the potential scaled to
$k=\tilde{k}\neq 1$ is given by
\begin{equation*}
\begin{cases}
\text{flock:} &  \tilde{\rho}(x) = \tilde{k}^{n}\rho(\tilde{k}x) \,,\, \operatorname{supp}(\tilde{\rho})=B(0,\frac{R}{\tilde{k}}) \\
\text{mill:} & \tilde{\rho}(x) = \tilde{k}^{2}\rho(\tilde{k}x) \,,\, \operatorname{supp}(\tilde{\rho})=B(\frac{R_{m}}{\tilde{k}},\frac{R_{M}}{\tilde{k}}) \\
\end{cases}.
\end{equation*}
Denoting $\tilde{U}(r) := U(\tilde{k}r)$, $\tilde W (x)=\tilde U
(|x|)$, we have $ \tilde{W}\star \tilde{\rho}= W\star \rho = C
\tilde{k}^{2-n}$ for flocks, and $\tilde{W}\star \tilde{\rho} =
W\star \rho + \frac{\alpha}{\lambda\,\beta}\log(\tilde{k})$ for
mill solutions.
\end{corollary}

\subsection{Explicit solvability }
In this section, we show how to solve almost explicitly the
integral equations for flock and mill profiles with the
Quasi-Morse potential. The exact problem to solve for any
potential $W$ is to find a density $\rho$ and its support such
that
\begin{equation}\label{charequationsec3}
(W \star \rho)(r) =  s(r) \text{ on } \operatorname{supp}(\rho)
\end{equation}
with some radial $s(r)$ on $\operatorname{supp}(\rho) =
B(R_m,R_M)$, $0\leq R_m < R_M$. Solving \eqref{charequationsec3}
generally implies inverting the integral operator, a task that is
complicated by the fact that the support is unknown. Even if this
is numerically achievable, we will not learn anything about the
structure of the solutions. For the Quasi-Morse potential, we take
advantage of the differential operators behind its construction,
to avoid the inversion of \eqref{charequationsec3} and to give an
almost explicit expression of its solution in terms of special
functions. This strategy was already done in \cite{BT2} in the one
dimensional case, where Morse and Quasi-Morse potentials coincide.
In this section, we pursue a similar strategy for the Quasi-Morse
potential in $n=2,3$.
\newline\indent
We begin our discussion reminding the radially symmetric
fundamental system associated with some operators, which will be
needed henceforth.
\newline\indent
\begin{remark} \label{remark-helmholtzsolutions}
The $n$-dimensional \emph{Helmholtz equation} reads $\Delta u +
k^{2} u = 0$ in $\rr^{n}$. Its fundamental system of radially
symmetric solutions $\{\varphi_1,\varphi_2\}$ is associated to a
second-order ordinary differential equation (in radial coordinates
for $n=2,3$) and given below, together with the fundamental system
of the already mentioned screened Poisson equations
$\{\psi_1,\psi_2\}$:
\begin{center}
\begin{tabular}{c c}
\begin{tabular}{c|c|c}
Helmh.& $\varphi_1$ & $\varphi_2$ \\
\hline $n=1$ & $\frac{1}{2k}\sin(kr)$ & $-\frac{1}{2k}\cos(kr)$ \\
\hline $n=2$  & $-\frac{1}{2\pi} J_{0}(kr)$ & $\frac{1}{2\pi}
Y_{0}(kr)$\\ \hline $n=3$ & $\frac{1}{4\pi} \frac{\sin(kr)}{r}$ &
$-\frac{1}{4\pi}\frac{\cos(kr)}{r}$
\end{tabular}
&
\begin{tabular}{c|c|c}
s.Poiss. & $\psi_1$ & $\psi_2$ \\
\hline $n=1$ & $\frac{1}{2k} e^{kr}$ & $- \frac{1}{2k} e^{-kr}$ \\
\hline $n=2$  & $\frac{1}{2\pi} I_{0}(kr)$ & $- \frac{1}{2\pi}
K_{0}(kr)$\\ \hline $n=3$ & $\frac{1}{4\pi} \frac{e^{kr}}{r}$ &
$-\frac{1}{4\pi}\frac{e^{-kr}}{r}$
\end{tabular}
\end{tabular}
\end{center}
Here, $J_0$ and $Y_0$ are the Bessel functions of the first and
second kind respectively, while $I_0$ and $K_0$ are the modified
Bessel functions of the first and second kind respectively.
\end{remark}
We continue with a simple computation related to the local
properties of our potential.
\newline\indent
\begin{lemma}
\label{lemma-laplacexoverl} Let $V$ be the fundamental solution
of the screened Poisson equation. Then
\begin{equation*}
\Delta_{x} \left(V\left(\frac{x}{l}\right)\right)=\frac{k^{2}}{l^{2}}V+l^{n-2}\delta_{0}.
\end{equation*}
\end{lemma}
\begin{proof}
Let $\xi$ be a test function. Then by change of variables
\begin{align*} \int_{\R^n}
\Delta_{x}\left(V\left(\frac{x}{l}\right)\right)\xi(x)\,\dd x &=
\frac{1}{l^{2}}\int_{\R^n} (\Delta
V)\left(\frac{x}{l}\right)\xi(x)\dd
x = \frac{1}{l^{2}}\int_{\R^n} \Delta V(z)\xi(lz)l^{n}\dd z \\
& =\frac{l^{n}}{l^{2}}\left(\int_{\R^n} k^{2}V(z)\xi(lz)\dd z +
\xi(0)\right) \\
&=\frac{1}{l^{2}}\int_{\R^n} k^{2}V(z)\xi(lz)l^{n}\dd z +
l^{n-2}\xi(0)
\\
&=\frac{k^{2}}{l^{2}}\int_{\R^n}
V\left(\frac{x}{l}\right)\xi(x)\dd x + l^{2}\xi(0)
\end{align*}
leading to the weak formulation of the claim.
\end{proof}
Now, we can state the main result of this section.
\newline\indent
\begin{theorem}\label{maincorollary}
Assume there exists a solution of $(W \star \rho)(r) =  s(r)$ on
$\operatorname{supp}(\rho)$ with $W$ being the Quasi-Morse
potential and $\operatorname{supp}(\rho)=B(0,R_{F})\,,s(r) = D$
for flocks, or $\operatorname{supp}(\rho)=B(R_{m,}R_{M})\,, s(r) =
D +\frac{\alpha}{\beta}\log(r)$ for mills respectively. Then
$\rho$ has to be of the following form on $\operatorname{supp}
\rho:$
\begin{center}
\begin{tabular}{c|c|c|c}
$n=2$: & flock & $A>0$ &$\rho_{F} = \mu_{1}\, J_{0}(ar) +\mu_{2}$\\
\hline
&  & $A=0$ &$\rho_{F} = \mu_{1}r^{2} +\mu_{2}$\\
\hline
& & $A<0$ &$\rho_{F} = \mu_{1}\, I_{0}(ar) +\mu_{2}$\\
\hline
& mill & $A>0$ & $\rho_{M} = \rho_{\text{inhom}}+ \mu_{1}\, J_{0}(ar) +\mu_{2}\, Y_{0}(ar)  +\mu_{3}$ \\
\hline
&  & $A=0$ & $\rho_{M} = \frac{\alpha}{\beta}\frac{k^{4}}{4\lambda l^{2}(1-C)}r^{2}(\log(r)-1) + \mu_{1}r^{2}  +\mu_{2}\log(r) +\mu_{3}$ \\
\hline
&  & $A<0 $& $\rho_{M} = \rho_{\text{inhom}}+ \mu_{1}\, I_{0}(-ar) +\mu_{2}\cdot K_{0}(ar)  +\mu_{3}$ \\
\hline
$n=3$: & flock & $A>0$ & $\rho_{F} = \mu_{1}\, \sin(ar)\frac{1}{r} +\mu_{2}$\\
\hline
 &  & $A=0$ &$\rho_{F} = \mu_{1}r^{2} +\mu_{2}$\\
\hline
 &  & $A<0$ & $\rho_{F}=\mu_{1}\sinh(ar)\frac{1}{r}+\mu_{2}$
 \\
\hline
\end{tabular}
\end{center}
with $A=k^{2}\frac{Cl^{n}-1}{l^{2}-Cl^{n}}$ , $a^2=|A|$, and
$\rho$ satisfying $\rho>0, \int\rho \dd x = 1$.
\end{theorem}
\begin{proof}
Let us define the operators $\mathcal{L}_{1}:=\Delta - k^{2}I$,
$\mathcal{L}_{2}:=\Delta -\frac{k^{2}}{l^{2}}I$. We apply both
operators to the equation and obtain
\begin{flalign*}
\mathcal{L}_{2}\mathcal{L}_{1}(W\star \rho)&= (\mathcal{L}_{2}\mathcal{L}_{1}W)\star \rho =\lambda  \left(-C\, \mathcal{L}_{1}\mathcal{L}_{2}V\left(\frac{r}{l}\right) +\mathcal{L}_{2}\mathcal{L}_{1}V(r) \right) \star \rho\\
&=\lambda\left( -Cl^{n-2}\Delta\delta+Ck^{2}l^{n-2}\delta +\Delta\delta -\frac{k^{2}}{l^{2}}\delta\right) \star \rho\\
&=\lambda(1-Cl^{n-2})\Delta \rho +\lambda\left(Ck^{2}l^{n-2}-\frac{k^{2}}{l^{2}}\right) \rho = \mathcal{L}_{2}\mathcal{L}_{1}s
\end{flalign*}
using Lemma \ref{lemma-laplacexoverl}. Hence, $\rho$ should
satisfy the following equation in its support:
\begin{equation}\label{resultinghelmholtz}
\Delta \rho \pm a^{2} \rho=\frac{1}{\lambda} \frac{1}{
1-Cl^{n-2}}\mathcal{L}_{2}\mathcal{L}_{1}s,
\end{equation}
with $a^2=|A|$ and
\begin{equation*}
A=\frac{Ck^{2}l^{n-2}-\frac{k^{2}}{l^{2}}}{1-Cl^{n-2}}=k^{2}\frac{Cl^{n}-1}{l^{2}-Cl^{n}},
\end{equation*}
resulting in the Helmholtz equation for $A>0$, the screened
Poisson equation for $A<0$ and the Poisson equation for $A=0$ with
radially symmetric inhomogeneous right-hand side. Therefore, the
solution to \eqref{resultinghelmholtz} writes as a general
solution of the homogeneous problem given by a linear combination
of the fundamental system in Remark
\ref{remark-helmholtzsolutions} plus a particular solution of the
inhomogeneous problem.
\newline\indent
The right-hand side of \eqref{resultinghelmholtz} depends on the
type of solution we wish to compute. For flocks in any dimension,
$s(r)$ is a constant function, and then we have
$\frac{1}{\lambda(1-Cl^{n-2})}\mathcal{L}_{2}\mathcal{L}_{1}s(r) =
\tilde{D}$. Therefore, the inhomogeneous solution of
\eqref{resultinghelmholtz} for $A\neq 0$ with unknown constant
right-hand side $\tilde{D}$ is
\begin{equation*}
\rho_{\text{inhom},A}(r)=\frac{\tilde{D}}{A}\,
1\hspace{-2pt}\text{\small I}_{\operatorname{supp} \rho}\,,
\end{equation*}
\newline\indent
For mills and $n=2$ we have $s(r)=D+\frac{\alpha}{\beta}\log(r)$
to obtain
\begin{flalign}
\frac{1}{\lambda(1-C)}\mathcal{L}_{2}\mathcal{L}_{1}\left[D+\frac{\alpha}{\beta}\log(r)\right]
&=\frac{k^{4}}{\lambda
l^{2}(1-C)}\frac{\alpha}{\beta}\log(r)+\tilde{D}
\label{righthandsidemill}
\end{flalign}
since $\log(r)$ is the fundamental solution of the Laplacian and
its Dirac delta disappears and we look for mill solutions on an
annulus (see Section \ref{section2}, \eqref{millcharrad}).
Therefore, the inhomogeneous solution of
\eqref{resultinghelmholtz} with right-hand side
\eqref{righthandsidemill} can also be written explicitly. Again
since $\log(r)$ is a fundamental solution of the Laplacian and the
support of the solution is assumed not to contain the origin, it
states
\begin{equation*}
\rho_{\text{inhom},A}(r) =  \frac{k^{4}}{\lambda
a^{2}l^{2}(1-C)}\frac{\alpha}{\beta}\log(r) +\frac{\tilde{D}}{A}
\text { on } \operatorname{supp} \rho \text{ for } A\neq 0.
\end{equation*}
Finally, in case $A=0$, the inhomogeneous solution for the flock
case is
\begin{equation*}
\rho_{\text{inhom},0}= \begin{cases}
\frac{1}{4}\tilde{D}r^{2} &, n=2 \\
\frac{1}{6}\tilde{D}r^{2} &, n=3
\end{cases} ,
\end{equation*}
whereas in the mill case it reads
\begin{equation*}
 \rho_{\text{inhom},0}=\frac{\alpha}{\beta}\frac{k^{4}}{4\lambda l^{2}(1-C)}r^{2}(\log(r)-1) + \frac{1}{4}\tilde{D}
 r^{2}\,
\end{equation*}
by using the fundamental solution of the Laplace operator. Putting
together the inhomogeneous solution with the homogenous part leads
to the claim of the theorem. For flocks, the space of candidate
solutions is of lower dimension, since singularities at the origin
are excluded.
\end{proof}
The coefficients $(\mu_{1},\mu_{2})$ or
$(\mu_{1},\mu_{2},\mu_{3})$ have to be computed numerically under
the constraint that the solution has to be non-negative, has to
contain unit mass, and has to solve the original equation
\eqref{charequationsec3}, but only on its own support which is a
priori unknown. To achieve this, we now need only to evaluate the
convolution integral in \eqref{charequationsec3} in a constrained
optimization method rather than its inversion.
\newline\indent
\begin{remark}\label{remflock3d}
Finally, we show that the radius of the support $R$ and the constants
$(\mu_{1},\mu_{2})$ are connected by an explicit nonlinear
identity in the particular case of 3D flocks. By plugging the
definition of the Quasi-Morse potential in 3D \eqref{potdef} into
\eqref{conv3d}, then
\begin{align*}
\Psi(r,s) = \frac{\lambda s}{2rk}\bigg( Cl^2
e^{-\frac{k}{l}|r-s|}-e^{-k|r-s|} - Cl^2
e^{-\frac{k}{l}(r+s)}+e^{-k(r+s)} \bigg)\,.
\end{align*}
By looking up in the table of Theorem \ref{maincorollary}, we have
the explicit expression of $\rho_F(r)$ for flocks with $A>0$ and
$n=3$. Straightforward computations lead to
\begin{align*}
\int_0^R \Psi(r,s)\rho_F(s) ds - \frac{\lambda \mu_2}{k^2}(Cl^3-1)
= &\,\frac{\lambda Cl^3}{k^3(k^2+a^2l^2)} \Lambda(C,l)
e^{-\frac{k}{l}R}\frac{\sinh\left(\frac{k}{l}r\right)}{r}\\
&-\frac{\lambda}{k^3(k^2+a^2)} \Lambda(1,1) e^{-kR}\frac{\sinh
kr}{r}
\end{align*}
where $A=a^2$ was used, and with
\[
\Lambda(C,l)=\mu_1k^2al\cos Ra+\mu_2
ka^2l^2R+l^3\mu_2a^2+\mu_2lk^2+\mu_1k^3\sin Ra + \mu_2k^3R.
\]
Therefore, the existence of a flock solution is equivalent to the
conditions $\Lambda(C,l)=0$ and $\Lambda(1,1)=0$, or equivalently
\[
 \begin{pmatrix}
k^2al\cos Ra + k^3\sin Ra & ka^2l^2R+a^2l^3+lk^2+k^3R \cr k^2a\cos
Ra + k^3\sin Ra &  ka^2R+a^2+k^2+k^3R
\end{pmatrix}
\begin{pmatrix}
 \mu_1 \cr
\mu_2
\end{pmatrix}
=\begin{pmatrix} 0 \cr 0\end{pmatrix}.
\]
A necessary condition ($\rho_F$ can still be negative) is that the
determinant of the matrix on the left hand side is zero, i.e., the
existence of the solution of the nonlinear equation for $R$
\[
 \tan Ra = \frac{a}{k} \frac{k^3R-a^2(l^2+l+klR)}{ka^2R+a^2(l^2+l+1)+k^2}.
\]
So that, for flock solutions in 3D we only need to check for radii
verifying this last identity.
\end{remark}
In the next section, we will show an algorithm to solve this
problem and present the numerical results.


\section{Numerical investigations}
\label{section4}

Theorem \ref{maincorollary} shows that solutions
of \eqref{charequationsec3} are solutions of
\eqref{resultinghelmholtz} with the constraints of positivity,
unit mass, and compact support on an annulus. Therefore, we now
propose an algorithm that numerically determines the support and
linear factors $\mu_{i}$ of the stationary flock and mill
solution. We will also present results which are compared to
particle simulations.

\subsection{The algorithm}
Let parameters $n,C,l,k,\alpha,\beta$ be fixed
and $\rho_{\text{hom}}$ denote the homogeneous solution dependent
on dimension as in Theorem \ref{maincorollary}. In search for the
support of the solution, we set the parameter $R_{max}$ as an
upper boundary on the support size the algorithms shall consider.
We can ensure that this is no restriction to the final result by
setting $R_{max}$ large compared to the characteristic shape of
the potential. Furthermore, denote $\Delta r$ a discretization
parameter and $\{r_{0},\ldots,r_{N}\}$ an equidistant
discretization of a chosen support $\text{supp}(\rho)$, with
$r_{i+1}-r_{i}=\Delta r$. Numerical approximations of functions
$F(r)$ on the discrete radial grid are denoted with $\bar{F}$. Our
first algorithm determines the best possible solution for
\emph{one particularly chosen} support $B(R_{l},R_{r})$. We aim to
find linear coefficients $(\mu_{1},\mu_{2})$ (or
$(\mu_{1},\mu_{2},\mu_{3}$) respectively), which solve the
integral equation \eqref{charequationsec3} the best possible way.
Non-negativity and unit mass of $\rho$ are hard constraints,
whereas the deviation $W\star \rho - s$ serves as the objective
function the coefficients shall minimize. 
\begin{algorithm}[for flocks]
\begin{flalign*}
\left|  \begin{array}{ll}
\textbf{Input}: \text{fixed support } B(0,R_{r})\\
\text{- For convolving functions with the potential, compute  a matrix } H  \text{ s.t. }\\ \,\,\, \overline{W\star \rho} = H\bar{\rho} \text{ according to Section \ref{sectionradialconvolutions}}.\\
\text{- Evaluate the convolution of the basis functions } \rho_{\text{hom}} \text{ and } \mathrm{1} \text{ on } \operatorname{supp} \rho: \\[2mm] \qquad\qquad\qquad g^{1}:=H\bar{\rho}_{\text{hom}}, g^{2}:=H\, \bar{\mathrm{1}}.  \\[3mm]
\text{-  To fit the right hand side $s(r)=D$ on the support, we chose coefficients} \\ \,\,\,\text{such that
$s(r)=D$ at the two end points $r_{1},r_{N}$. That is, solving }\\[2mm]
\qquad\qquad\qquad \left.\begin{pmatrix}
g^{1}_{1} & g^{2}_{1} \\ g^{1}_{N} & g^{2}_{N}
\end{pmatrix}\right. \mu_{\text{const}}= \begin{pmatrix}
1 \\ 1
\end{pmatrix}\\[3mm]
\,\,\, \text{setting } D=1 \text{ temporarily}.\\
\text{- By linearity of $H$, we set } \bar{\rho}:= \frac{1}{M}(\mu_{\text{const},1}\rho_{\text{hom}}+\mu_{\text{const},2}) \text{ with } M \text{ normalizing total mass.}
\\
\text{- Since we have only ensured \eqref{charequationsec3} to hold at two points, we measure deviation of } \\ \,\,\,\text{$H\bar{\rho}$ from $s(r)$ (here, an arbitrary constant) on the whole support as }\\[2mm]
\qquad\qquad\qquad \displaystyle{e:=\frac{1}{R_{r}}\int \left[H\bar{\rho} - \frac{1}{R_{r}}\int H\bar{\rho}\dd \bar{r} \right] \dd \bar{r}}.\\[3mm]
 \textbf{Output}: e, \bar{\rho}, \bar{s} \text{ if } \bar{\rho} \geq 0, \text{error message if } \bar{\rho} \ngeq 0.
\end{array} \right.
\end{flalign*}
\label{algo1}
\end{algorithm}

 For the case of mills, we proceed analogously, but we have to take into account the fixed inhomogeneous solution and three basis functions.
\begin{algorithm2}[for mills]
\begin{flalign*}
\left|  \begin{array}{ll}
\textbf{Input}: \text{fixed support } B(R_{m},R_{M})  \\
\text{- For convolving functions with the potential, compute  a matrix } H \text{ s.t. }\\
\,\,\,\overline{W\star \rho} = H\bar{\rho} \text{ according to Section \ref{sectionradialconvolutions}}.\\
\text{- Evaluate the convolution of the fixed inhomogeneous part } \rho_{\text{inhom},A} \text { on } \operatorname{supp} \rho \\
\,\,\,\text{ and set } \bar{s}_{\text{inhom}}:=H\bar{\rho}_{\text{inhom},A}.\\
\text{- Define the remainder of the right-hand side as } \bar{s}_{\text{rem}}:=\bar{s}-\bar{s}_{\text{inhom}},\\
\,\,\, \text{which has to be fitted by the convolution of the basis functions.}\\
\text{- To do so, evaluate } J_{0}(ar), Y_{0}(ar)\text{ and } \mathrm{1} \text{ on } \operatorname{supp} \rho \text{:}\\[2mm]
\qquad\qquad\qquad  g^{1}:=H\bar{J_{0}}, g^{2}:=H\bar{Y_{0}} ,g^{3}:=H \bar{\mathrm{1}}.  \\
\text{- Giving three basis functions, we pick three points $r_1, r_j \text{ with }j={\lfloor N/2\rfloor},r_{N}$ } \\
\,\,\,\text{and interpolate both the remainder $\bar{s}_{\text{rem}}$ and the free constant, which is }\\
\,\,\,\text{temporarily set to $1$. We solve }   \\[2mm]
\qquad\qquad\qquad\vspace{1mm} \left.\begin{pmatrix}
g^{1}_{1} & g^{2}_{1} & g^{3}_{1} \\ g^{1}_{j} & g^{2}_{j} & g^{3}_{j}  \\
g^{1}_{N} & g^{2}_{N} & g^{3}_{N}
\end{pmatrix}\right.\mu_{\text{rem}}= \begin{pmatrix}
\bar{s}_{\text{rem},1} \\ \bar{s}_{\text{rem},j} \\\bar{s}_{\text{rem},N}
\end{pmatrix}  \\ \vspace{1mm}
\qquad\qquad\qquad \left.\begin{pmatrix}
g^{1}_{1} & g^{2}_{1} & g^{3}_{1} \\ g^{1}_{j} & g^{2}_{j} & g^{3}_{j} \\
g^{1}_{N} & g^{2}_{N} & g^{3}_{N}
\end{pmatrix}\right.\mu_{\text{const}}= \begin{pmatrix}
1 \\ 1 \\1
\end{pmatrix} \\[4mm]
\text{- By linearity of $H$, we set } \bar{\rho}_{\text{rem}}:=\mu_{\text{rem},1}\bar{J_{0}}+\mu_{\text{rem},2}\bar{Y_{0}}+\mu_{\text{rem},3} \text{ and }  \\[2mm] \,\,\,\, \qquad\qquad \bar{\rho}_{\text{const}}:=\mu_{\text{const},1}\bar{J_{0}}+\mu_{\text{const},2}\bar{Y_{0}}+\mu_{\text{const},3}.\\[2mm]
\text{- Our last degree of freedom is the free constant on the right hand side $s(r)$, }\\
\,\,\,\text{which we use to normalise mass. The candidate density is} \\[2mm]
\qquad\qquad \bar{\rho} := \bar{\rho}_{\text{inhom},A}+\bar{\rho}_{\text{rem}}+\gamma\bar{\rho}_{\text{const}}\text{ with }\gamma := \frac{1-m(\bar{\rho}_{\text{rem}})-m(\bar{\rho}_{\text{inhom},A})}{m(\bar{\rho}_{\text{const}})} .\\[2mm]
\text{- We penalize deviation of $H\bar{\rho}$ from $s(r)$ on the entire support as} \\
\quad\quad\quad \displaystyle{e_{1}:=\frac{1}{R_{M}-R_{m}}\int \left[H\bar{\rho} - \bar{s} - \frac{1}{R_{M}-R_{m}}\int (H\bar{\rho} - \bar{s}) \dd \bar{r} \right] \dd \bar{r}}.\\
 \text{- Second, since $s(r)$ is concave, we penalize numerical convexity of } \bar{s} \text{ by }\\
 \qquad\qquad\qquad \displaystyle{e_{2}:= \int \chi_{[\bar{s}''>0]} \, \bar{s} \dd\bar{r}} \\
 \,\,\, \text{The total penalty value is the sum of $e_1, e_2$.}\\
 \textbf{Output}: e=e_{1}+e_{2}, \bar{\rho}, \bar{s} \text{ if } \bar{\rho} \geq 0
\end{array} \right.\\
\end{flalign*}
\end{algorithm2}
Now, we search the minimizer of
the error function $e$ over a test set of supports, given by
the pre-defined discretization $\Delta r_1$and maximal support size. Repeating Algorithm \ref{algo1} over the set of test supports provides a minimizer of the penalty function. For flocks, the number of tested supports is  $\approx \frac{R_{\text{max}}}{\Delta r_1}$, for mills  $\approx \frac{1}{2}\left(\frac{R_{\text{max}}}{\Delta r_1}\right)^{2}$. To enhance the speed of numerical computation, we first compute a solution based on a coarser discretization length $\Delta r_{2}=m \Delta_{r_1}$ for some integer $m$. Then, the obtained minimizer is used as the center of a local refinement search with a fine discretization length, as illustrated in Algorithm \ref{algo2}:
\begin{algorithm}
\begin{flalign*}
\left|  \begin{array}{ll}
\text{- Choose a coarse grid size }\Delta r_{2} \text{ such that an iteration of Algorithm \ref{algo1} over all}\\\text{test supports is reasonably fast and, as a solution, obtain the support } B(0,\tilde{R}_{F}) \\
\,\,\, \text{(or } B(\tilde{R}_{m}, \tilde{R}_{M})) \text{ for mills}. \\
\text{- Vary this support locally} \text{ up to a fixed parameter $c$ with a fine discretization }\\
\,\,\, \Delta r_{1}\ll \Delta r_{2} \text{, re-run Algorithm \ref{algo1} restricted on } |{R}_{F}-\tilde{R}_{F}|\leq c \\
\text{\,\,\,\,} (|R_{m }-\tilde{R}_{m}|\leq c, |R_{M }-\tilde{R}_{M}|\leq c   \text{ for mills).}\\
\text{- Obtain the minimising results  supp$(\rho), \bar{\rho}, \bar{s}$ and $e$. }
\end{array} \right.
\end{flalign*}
\label{algo2}
\end{algorithm}

Naturally, the matrix $H$ is not recomputed in every iteration but
constructed once for the largest support and inherited. The choice
to fix a functional equality on the points which are most left,
most right and for mills central on the chosen support is
arbitrary. We say that no compact solutions are found in our
computations, if our algorithms deliver $R_{\text{max}}$ as the
error minimizer, no matter of its value. The convergence of the
algorithm for $\Delta r \rightarrow 0$ if compact solutions are
found will be demonstrated together with the results of the next
subsection.


\subsection{Flocks in 2D}
\label{sec-subflocks2d}

\begin{figure}
\subfloat[Flock emerged in a particle simulation with  $N=400$
particles]{\includegraphics[keepaspectratio=true,width=.5\textwidth]{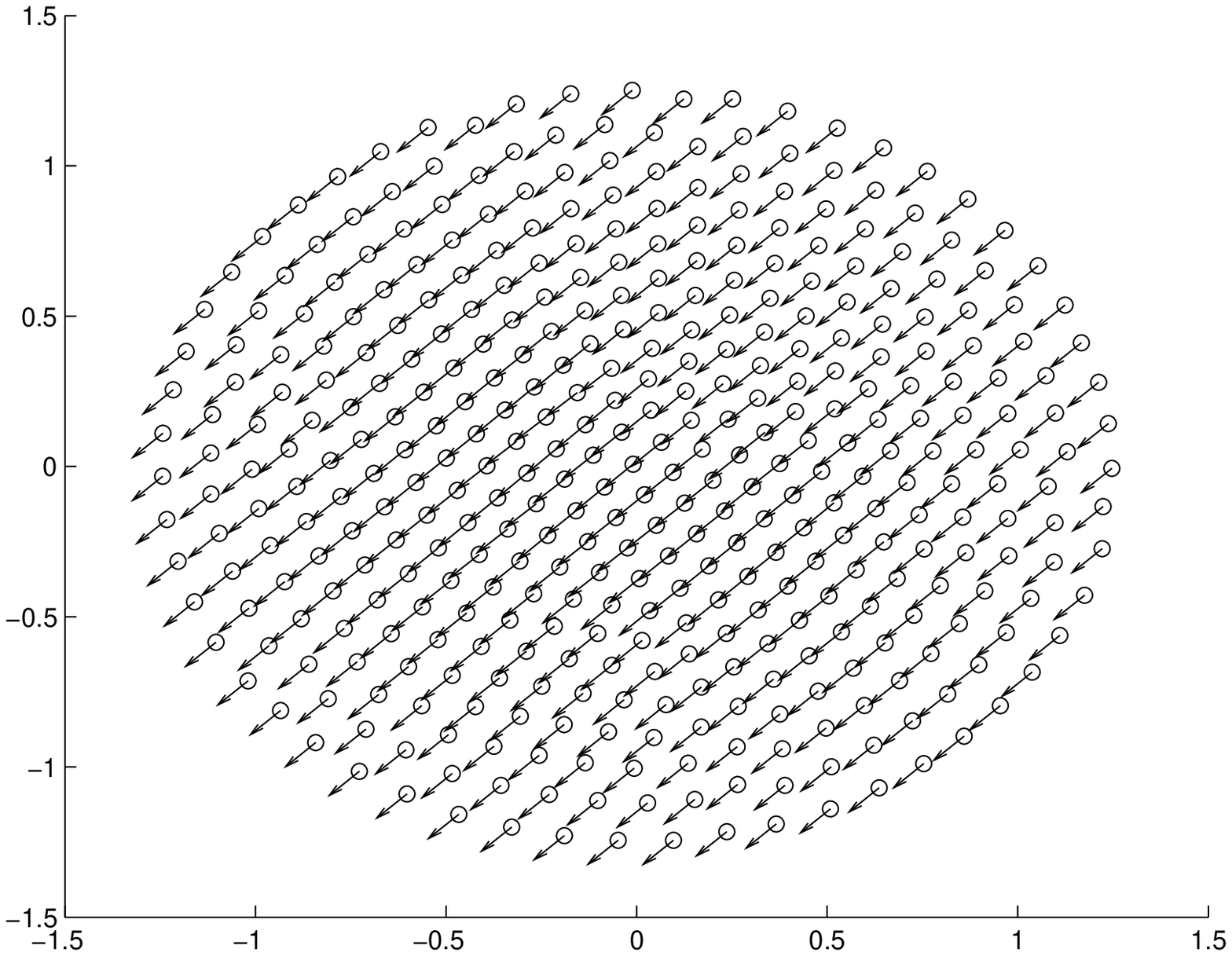}}
\hspace{.05cm} \subfloat[Radial flock density: Continuous result
vs. empirical density ($N=30000$
particles)]{\includegraphics[keepaspectratio=true,width=.5\textwidth]{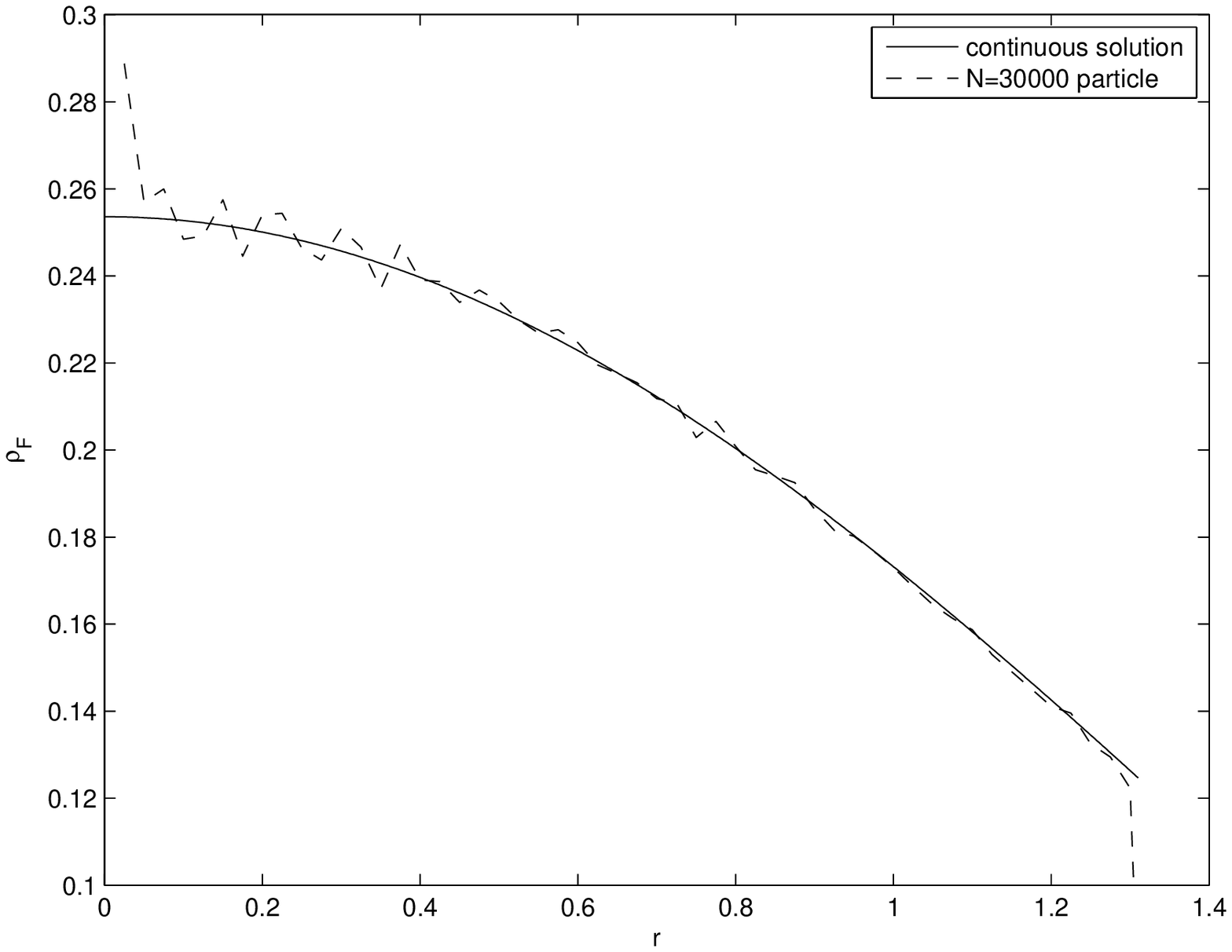}}
\caption{Two-dimensional aligned flocks emerge for the Quasi-Morse
potential. The resulting continuous radial density of Algorithms
\ref{algo1},\ref{algo2} matches the empirical distribution
obtained from particle simulations. The stationary flock has the
form $\rho_{F} = \mu_{1}\, J_{0}(ar) +\mu_{2}$ with, in this case,
$\mu_{1}\approx 0.2356, \mu_{2}\approx0.018, A=1.5,
R_{F}\approx1.31$ (Quasi-Morse potential parameters in use are
$C=\frac{10}{9}, l=0.75, k=\frac{1}{2}$).} \label{fig-flock2d}
\end{figure}

\begin{figure}
\subfloat[Continuous solution $\rho_{F}$ for varying $\Delta r$]{\includegraphics[keepaspectratio=true,width=.5\textwidth]{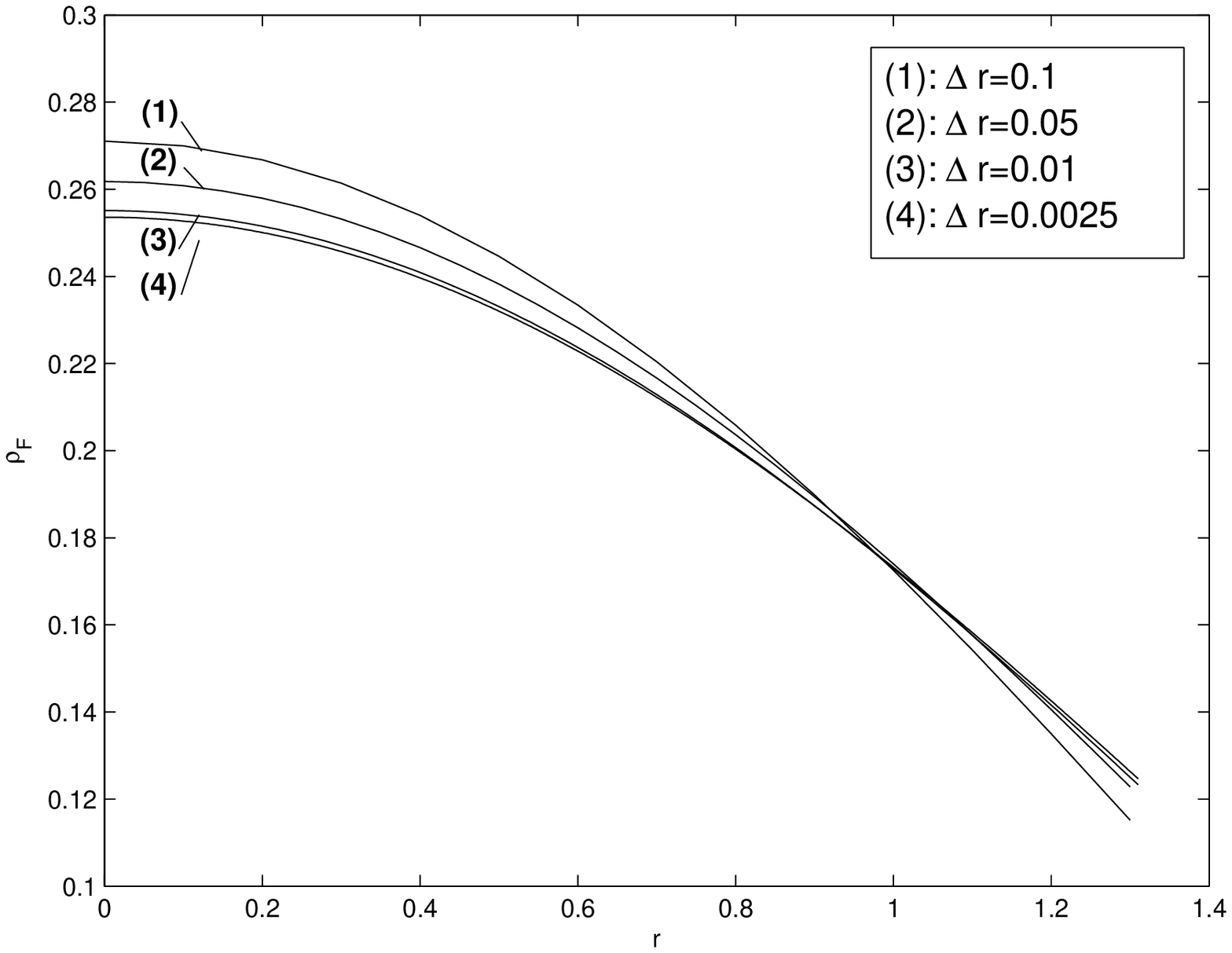}}
\hspace{.05cm}
\subfloat[Minimal error value and computation times]{
\centering
\vspace{-2cm}
\begin{tabular}[b]{|c|c|c|} \hline $\Delta r $ & error $e$ & computation time  \\ \hline
0.1 & 3.54e-05  & 0.76s \\
0.05 & 1.36e-05 & 2.85s \\
0.01 & 3.99e-06 & 69.1s \\
0.0025 & 9.97e-07 & 1125s \\
\hline
\end{tabular}}
\caption{Algorithms \ref{algo1},\ref{algo2} converge as $\Delta r
\rightarrow 0$ if a compactly supported flock solution exists.
Four resulting densities are shown for $\Delta r
\in\{0.1,0.05,0.01,0.0025\}$ together with the minimal error value
of the algorithm and the corresponding computation time.}
\label{fig-flock2dconvergenceandtable}
\end{figure}

We start our presentation of numerical results with the aligned
flock in two dimensions. Our standard example is the configuration
$C=\frac{10}{9}, l=0.75, k=\frac{1}{2}$ as in Fig.
\ref{fig-comparepotentials}. The stationary aligned flock state is
independent of $\lambda, \alpha, \beta$, yet emergence of flocks
in particle simulations depends on these parameters and suitable
initial conditions. An exemplary convenient choice is $\alpha=1,
\beta=5, \lambda\in\{100,1000\}$. The observed flock of aligned
particles is illustrated in Fig \ref{fig-flock2d}a for $N=400$
particles. In Fig. \ref{fig-flock2d}b, the result of our
investigations is compared to the empirical radial density
obtained from a particle simulation with $N=30000$ agents.
The empirical radial density is obtained by
collecting particles in radial bins and dividing by the Jacobian
of the radial transformation. We see that the continuous solution
matches the particle density and convergence is expected as
$N\rightarrow\infty$.
While the numerical cost of full particle simulations is at least
$\mathcal{O}(N^{2})$, the computational effort of the presented
method scales quadratically with $\Delta r$, as illustrated in
Fig. \ref{fig-flock2dconvergenceandtable}b. In Fig.
\ref{fig-flock2dconvergenceandtable}a we show the convergence of
our algorithms as $\Delta r \rightarrow 0$. One observes that the
support is estimated well for coarse grid sizes, whereas the
correct radial density is established with finer discretizations.
The minimal error values of Algorithms \ref{algo1}, \ref{algo2}
are listed in Fig. \ref{fig-flock2dconvergenceandtable}b.
The advantages of the presented solution are continuity, dramatic
reduction of the numerical cost, fast convergence, and an explicit
expression of the radial density as, in this example, a
combination of Bessel's J-function and a constant.

\begin{figure}
\centering
\includegraphics[keepaspectratio=true, width=.5\textwidth]{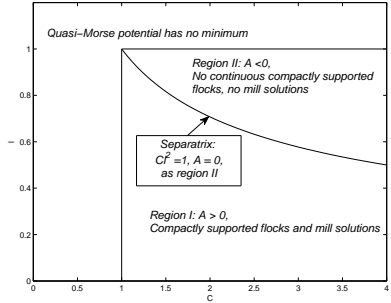}
\caption{Numerical phase diagram of the
Quasi-Morse potential in 2D: The biologically relevant scenarios
decompose into two subregions. Region I: $A>0$, continuous
compactly supported flocks. Region II: $A\leq 0$, no compactly
supported continuous solutions, flocks only emerge on particle
level. The same division of regions applies to the mill
solutions.} \label{fig-2dflockparameters}
\end{figure}

Concerning the potential parameters, the area of relevant
short-term repulsion and long-range attraction shapes divides into
two subregions based on the results of section \ref{section3}, as
illustrated in Fig. \ref{fig-2dflockparameters}: In region I with
$C>1,l<1, Cl^{2}<1$, the potential is catastrophic, $A>0$
(from Theorem \ref{maincorollary}) and compactly
supported continuous flock solutions are found. In region II with
$C>1,l<1, Cl^{2}<1, A<0$ (from Theorem
\ref{maincorollary}), and no compactly supported solutions can be
found. No solution of the algorithm indicates
H-stability since in this case particle simulations do show flocks
whose support diverges when $N\to\infty$. The presented method
faces numerical difficulties for catastrophic potentials
$A>0$ with $Cl^{2}\approx 1$, where it eventually
breaks down not converging to a compactly supported flock.
Similarly, particle simulations are not fully reliable in this
limiting cases. However, thanks to our computation in Section
\ref{section3} we are able to consider the exact
separatrix case $Cl^{2}=1, C>1, l<1, A=0$: Here, no compact
solutions are found. Our numerical findings are
illustrated in Fig. \ref{fig-2dflockparameters}. We emphasize that
based on the reported simulations, we conjecture that compactly
supported flock solutions exist only in the catastrophic regime,
$A>0$.

\subsection{Mills in 2D}

\begin{figure}
\subfloat[Mill emerged in a particle simulation with  $N=400$
particles]{\includegraphics[keepaspectratio=true,width=.5\textwidth]{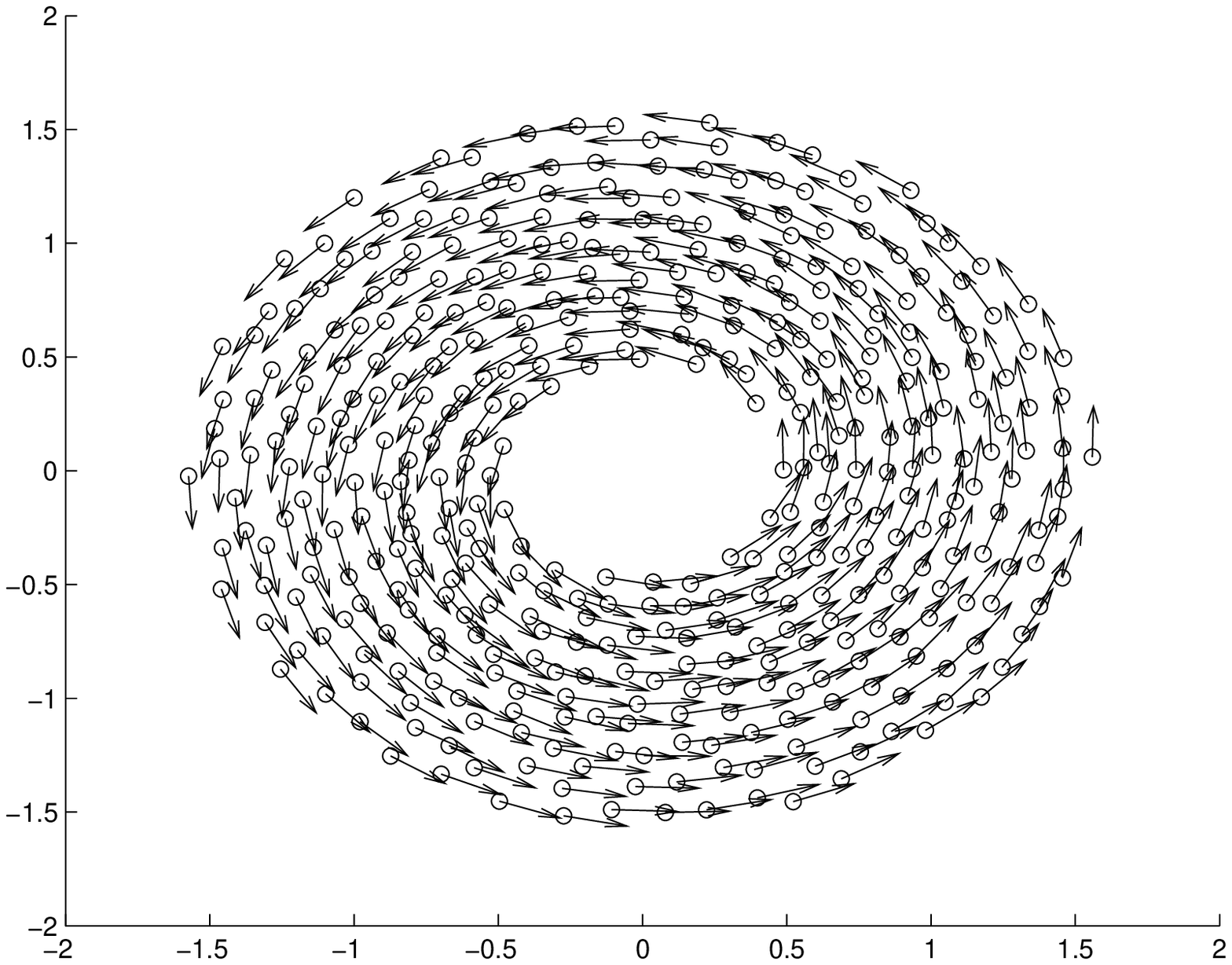}}
\hspace{.05cm} \subfloat[Radial mill density: Continuous result
vs. empirical measure ($N=16000$
particles)]{\includegraphics[keepaspectratio=true,width=.5\textwidth]{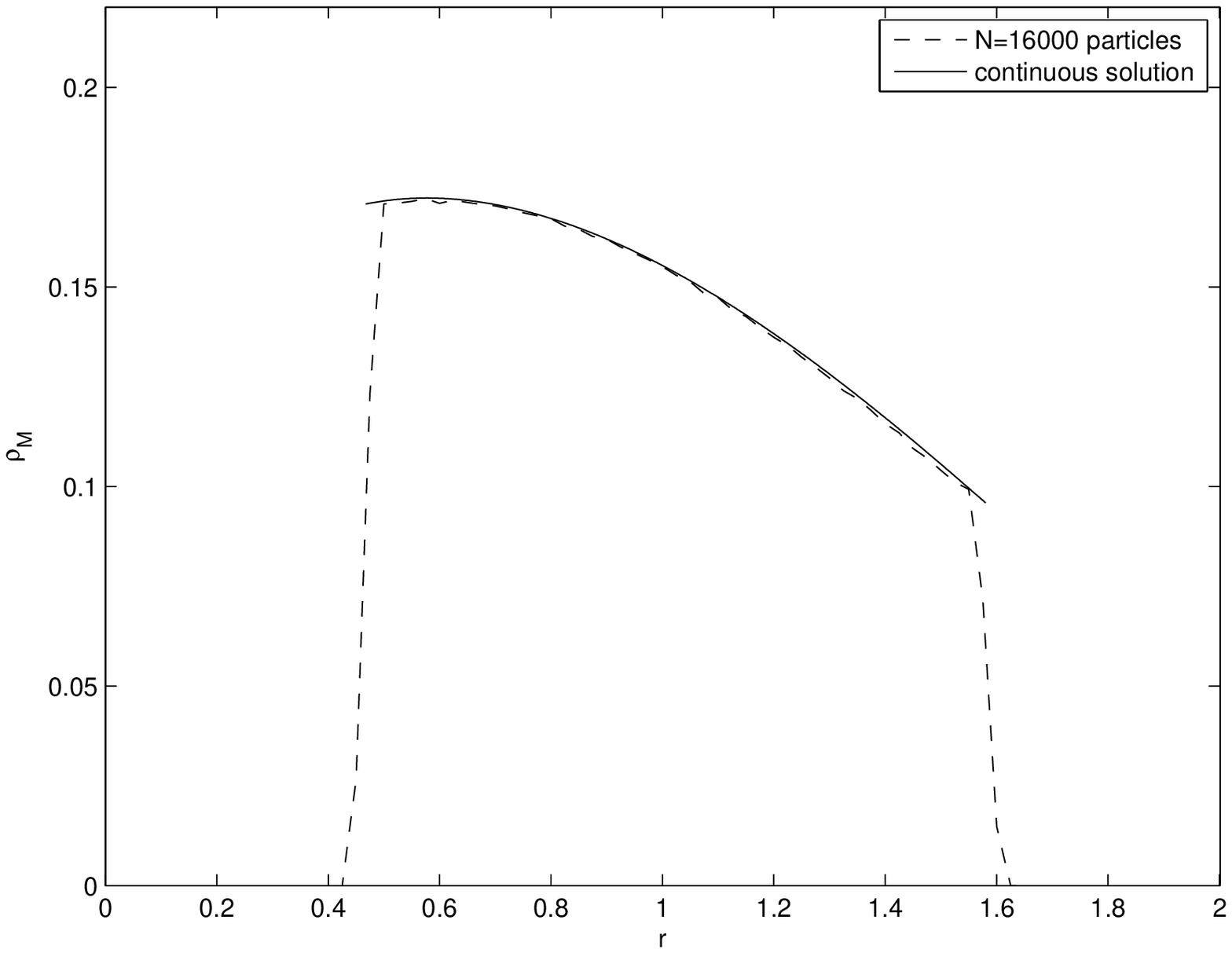}}
\caption{Rotating mills emerge for the Quasi-Morse potential. As
for flocks, the resulting radial density of Algorithms
\ref{algo1},\ref{algo2} matches the empirical distribution
obtained from particle simulations. The  mill solution has the
form $\rho_{M} = \rho_{\text{inhom},A}+ \mu_{1}\, J_{0}(ar)
+\mu_{2}\, Y_{0}(ar)  +\mu_{3}$ with, in this case,
$\mu_{1}\approx 0.1708, \mu_{2}\approx 0.0468, \mu_{3}=0.0320,
A=1.5, \text{supp}_{\rho_{M}}\approx B(0.47, 1.57)$ (Quasi-Morse
potential parameters in use are $C=\frac{10}{9}, l=0.75,
k=\frac{1}{2}$, others are $\alpha=1, \beta=5, \lambda=100$).}
\label{fig-mill2d}
\end{figure}

\begin{figure}
\centering
\includegraphics[keepaspectratio=true,width=.5\textwidth]{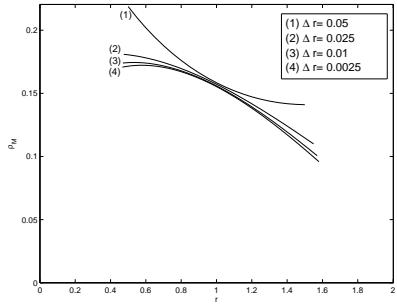}
\caption{Algorithms \ref{algo1},\ref{algo2} converge for the mill
case as $\Delta r \rightarrow 0$. } \label{fig-mill2dconv}
\end{figure}

\begin{figure}
\subfloat[mill solutions for varying
$\frac{\beta\lambda}{\alpha}$]{\includegraphics[keepaspectratio=true,width=.5\textwidth]{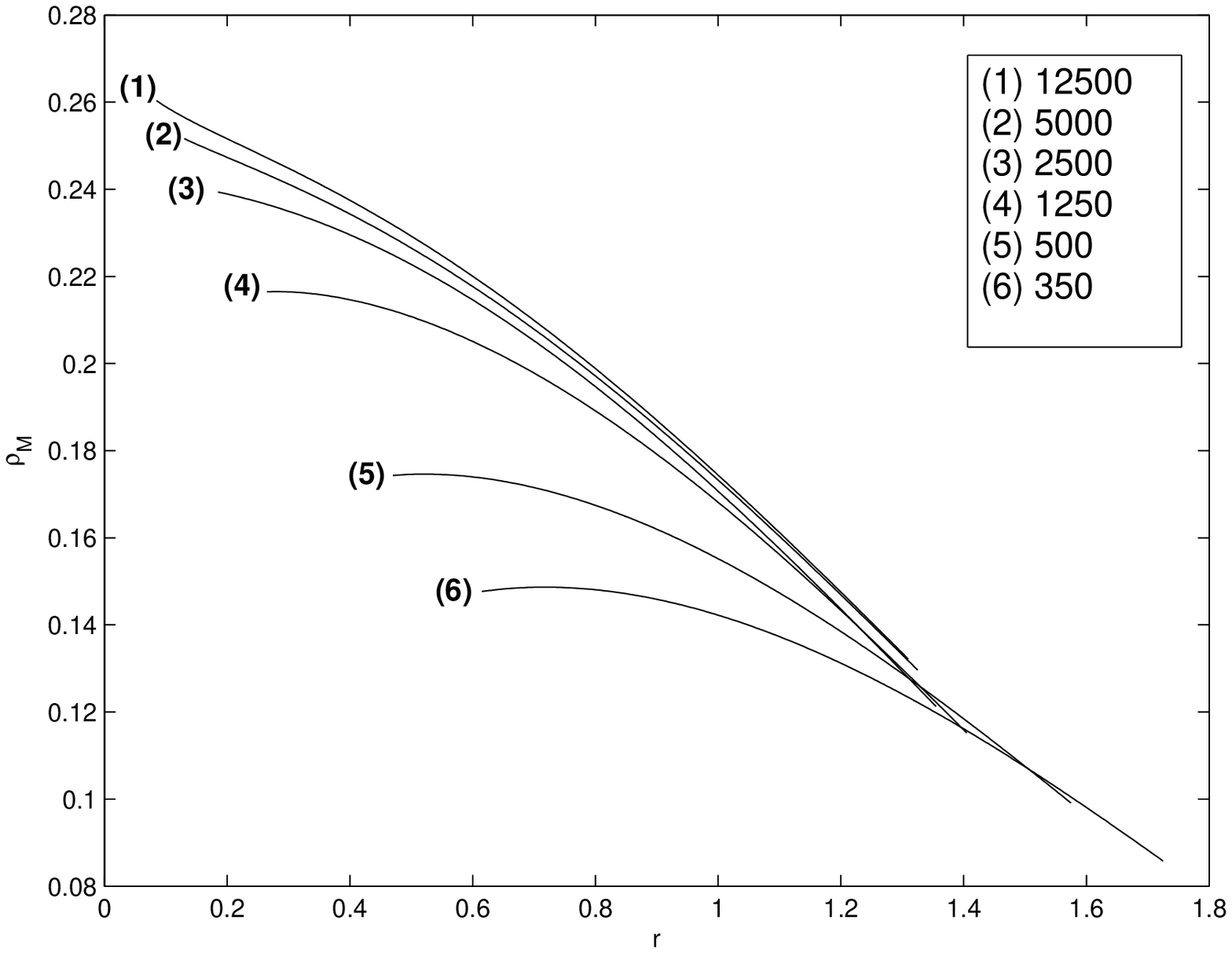}}
\hspace{.05cm} \subfloat[support of mill solutions for varying
$\frac{\beta\lambda}{\alpha}$]{\includegraphics[keepaspectratio=true,width=.5\textwidth]{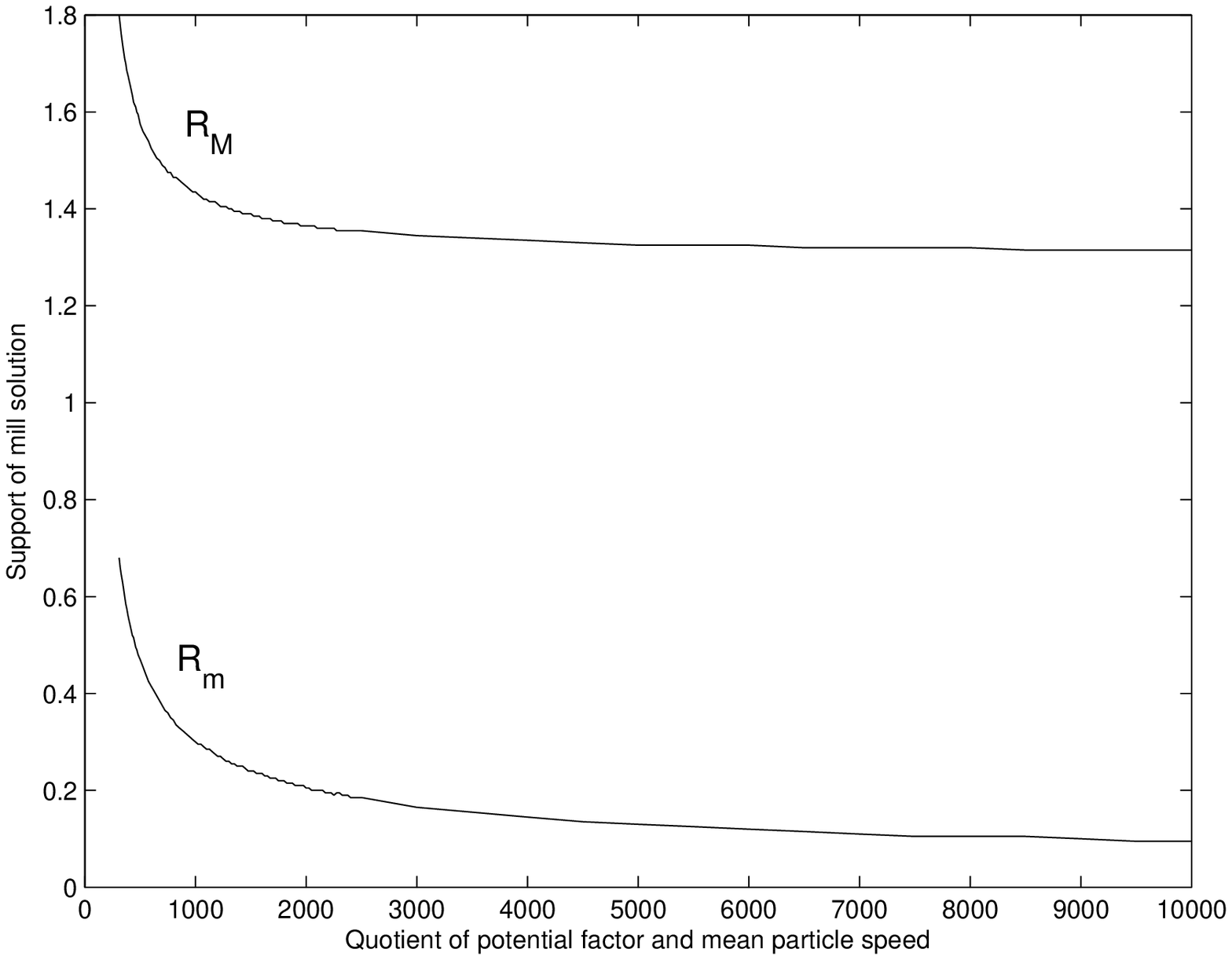}}
\caption{Quasi-Morse potentials with identical shape parameters
$C,l,k$ result in mill solutions with different support sizes and
densities, depending on the ratio of potential factor and squared
stationary speed of the mill. }
\label{fig-mill2d-quotientinfluence}
\end{figure}

The Quasi-Morse potential is able to produce rotating mill states
in particle simulations, just as the original Morse potential. We
choose the same configuration as in Section \ref{sec-subflocks2d}
with $\lambda=100$ and show the mill emerging from a particle
simulation in Figure \ref{fig-mill2d}a. The resulting mill
solution of our algorithms is illustrated in Figure
\ref{fig-flock2d}b, together with a comparison to an empirical
density from a particle mill with $N=16000$ agents. Again, our
result is confirmed by the particle simulation and support as well
as the density shape agree perfectly. The stationary rotating mill is
a weighted sum of Bessel's J and Y functions, the inhomogeneity
$\rho_{\text{hom}}$ and a constant.
The convergence of Algorithms \ref{algo1}, \ref{algo2} in the mill
case is shown in Figure \ref{fig-mill2dconv}. As for flocks, the
computational costs are minimal compared to a full particle
simulation.
For the existence of compactly supported mill solutions, the
parameter diagram on Figure \ref{fig-2dflockparameters} applies
just as  for flocks. In region I, continuous solutions can be
found, whereas in 
region II and the separatrix
$Cl^{2}=1$ no such mills can be found. In particle simulations, we
there see either a crystal-like arrangements or ''finite
particle'' flocks as in Section \ref{sec-subflocks2d}. Next we
study the impact of parameters $\alpha,\beta,\lambda$ on the
stationary mill solution, which enter the solution solely in the
joint quotient $\frac{\alpha}{\lambda\, \beta}$. Hence, for a
potential multiplied by a factor $\lambda$, the mill solution will
stay the same, if the preferred speed of particles is multiplied
by $\sqrt{\lambda}$ by any suitable change of $\alpha$ and/or
$\beta$. In Figure \ref{fig-mill2d-quotientinfluence}a, we show
several mill densities for our standard potential configuration
and $\frac{\beta}{\lambda \, \alpha}\in \{350, 500, 1250, 2500,
5000, 12500\}$. The support of mill solutions is plotted against
$\frac{\beta}{\lambda \, \alpha}$ in Figure
\ref{fig-mill2d-quotientinfluence}b.

\subsection{Flocks in 3D}
\begin{figure}
 \subfloat[Quasi-Morse potential in 3D]{
\includegraphics[keepaspectratio=true, width=.5\textwidth]{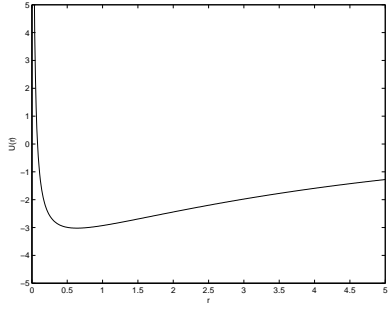}
}
 \hspace{.05cm}
 \subfloat[3D flock emerged in a particle simulation with $N=200$ particles]{
\includegraphics[keepaspectratio=true, width=.5\textwidth]{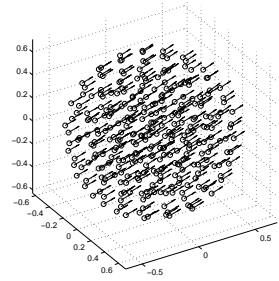}
}
\\
\subfloat[3D radial flock density: continuous result vs. empirical measure ($N=35000$ particles)]{
\includegraphics[keepaspectratio=true, width=.5\textwidth]{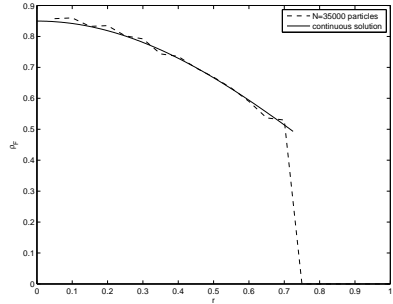}
}
 \hspace{.05cm}
 \subfloat[Parameter diagram]{
\includegraphics[keepaspectratio=true, width=.5\textwidth]{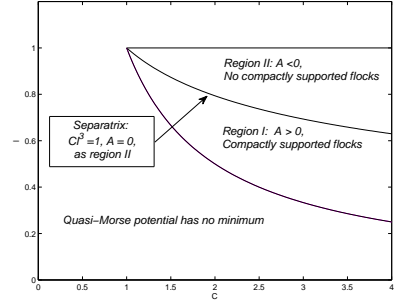}
} \caption{The Quasi-Morse potential in three dimension is able to
produce aligned flock solutions. The continuous radial density can
be expressed as $\rho_{F} = \mu_{1}\cdot \sin(ar)\frac{1}{r}
+\mu_{2}\cdot\mathrm{1}$ with, in this case, $\mu_{1}\approx
0.3574, \mu_{2}\approx 0.0052, R_{F}\approx 0.725, A=5.585$. Our
result is verified by comparing to the empirical density obtained
from a particle simulation.   (a) Exemplary potential shape, (b)
Flock emerged from 3D particle simulation, (c) Continuous flock
solution vs empirical measure, (d) Parameter diagram of biological
relevant configurations} \label{fig-3d}
\end{figure}

The introduction of Quasi-Morse potentials enables us also to
study flocks in three space dimensions. As we have mentioned in
Section \ref{section3}, the area of admissible parameter
configurations is smaller than in the 2D case, as illustrated in
the parameter diagram. For our example, we set $C=1.255, l=0.8,
k=0.2, A=5.585$ and plot the resulting potential shape in Figure
\ref{fig-3d}a. A three-dimensional flock resulting from a particle
simulation is shown in Figure \ref{fig-3d}b. With the help of
Algorithm \ref{algo1} the continuous radial flock density is
computed as a linear combination of $\frac{\sin ar}{r}$ and a
constant. Notice that due to Remark
\ref{remflock3d}, Algorithm \ref{algo2} is not needed. Also in
three dimensions, the empirical density of a particle simulation
matches our result, as illustrated in Figure \ref{fig-3d}c.
Concerning the existence of flock solutions in dependence of the
shape parameters $C$ and $l$, we get an equivalent picture as in
two dimensions (see Figure \ref{fig-3d}d): Though different in
shape, the area of biologically relevant shapes is divided into
two subregions by the separatrix $Cl^{3}=1$. In region I,
continuous compactly supported three-dimensional flocks are found,
not so in region II, which again indicates H-stability. Here,
flocks do appear but their support increases with the total number
of agents $N$. In the special case of the separatrix, which can be
investigated with the computation of the case $A=0$ in Section
\ref{section3}, no flock solutions are found.

\section{Discussion}
Quasi-Morse potentials fulfill three properties desirable from
biological modeling: short-term repulsion, long-term attraction
and vanishing interaction at infinity. Using Quasi-Morse
potentials instead of the standard Morse potential makes, in our
view, hardly any difference in terms of biological modeling. The
stronger singularity at the origin for $n\neq 2$ might even be
desirable in order to enforce repulsion. Though the special
functions involved for $n=2$ may seem not as convenient to work
with as the exponential function, the existence of continuous,
compactly supported stationary states itself make Quasi-Morse
potentials a good choice for further studies of the models
discussed in the above. Our results are, to the best of our
knowledge, one of the first of its kind for explicit solutions of
flock and mill patterns in two or three dimensions. The strategy
of building up potentials from solutions of certain partial
differential equations might work in other cases as well and form
one tool in the effort to understand the equilibria of interaction
potentials. However, the techniques applied here are of no help
for general potentials, such as classical Morse. With a variety of
potentials suggested (see the discussion in Section
\ref{section2}), the problem of choosing the best suited one for a
particular biological application becomes increasingly evident and
should be a topic of future research.

\section[nonumber]{Conclusions}
In this paper we have introduced the Quasi-Morse interaction
potentials for a second-order model of self-propelled interactive
particles. The Quasi-Morse potentials lead to the emergence of
flocks and mills, similar to the standard Morse potential.
We have shown that the radial densities
of these stationary states are (affine) linear combinations of two
or three elementary functions, which are chosen with the respect
to the three subcases $A>0$ (catastrophic), $A=0$ (separatrix) or
$A<0$. In order to determine the correct scalar
coefficients and the a priori unknown support, we have developed a
numerical algorithm that does not use time evolutions in Section
\ref{section4}. We have illustrated our result with examples for
flocks and mills in two dimension and flocks in 3D. In all cases,
our findings are convincingly verified by corresponding particle
simulations. With our algorithm, we find that for all coherent
patterns, only the catastrophic scenarios $A>0$ lead to continuous
compactly supported solutions.


\section{Appendix}
This appendix is devoted to the proof of Corollary
\ref{cor-biorel}. Denote by $V_k(r)$ the potential in
\eqref{potdef} for a given $k>0$. Notice that $V_k(r)=k^{n-2}
V_1(kr)$, hence we set $k=1$ without loss of generality and, from
now on, we drop the index $k$ for simplicity.
\newline\indent
Let us first show the assertion about the unique minimum of the
potential. This property is desirable from the biological point of
view to set a typical length scale for the distance between
agents.
\newline\indent
Let \(U(r)=V(r)-CV\big(\frac{r}{l}\big)\), then
\(U'(r)=V'(r)-\frac{C}{l} V'\big(\frac{r}{l}\big)\), and the
necessary condition for a local extremum can be stated as finding
$r$ such that
\[
h(r):=\frac{C}{l}\,\frac{V'\big(\frac{r}{l}\big)}{V'(r)} =1.
\]
Straightforward computations lead to
\[
\begin{split}
h'(r) & =\frac{C}{l(V'(r))^{2}}\big( \tfrac{1}{l}
V''(\tfrac{r}{l})V'(r) - V'(\tfrac{r}{l})V''(r)\big)
\\
& = \frac{C}{l(V'(r))^{2}}\big(
-\tfrac{n-1}{r}V'(r)V'(\tfrac{r}{l})+\tfrac{1}{l}V(\tfrac{r}{l})V'(r)+\tfrac{n-1}{r}V'(r)V'(\tfrac{r}{l})
\\
& -V(r)V'(\tfrac{r}{l})\big) = \frac{C}{l(V'(r))^{2}}\big(
\tfrac{1}{l}V'(r)V(\tfrac{r}{l})-V'(\tfrac{r}{l})V(r)\big),
\end{split}
\]
where we used that \(V\) is a solution of
\(\frac{n-1}{r}V'(r)+V''(r) = V(r)\).
\newline\indent
We next check that \(\log(-V(r))\) is a convex function of \(r\)
for \(n=1,2,3\). Indeed, if \(n=1\), then \(\log(-V(r))\) is
affine; if \(n=3\), then \(\log(-V(r))=-\log r - r -\log(4\pi)\).
In the case \(n=2\), we have
\[
(\log(-V(r)))'' =
\frac{K_0(x)^2+K_2(x)K_0(x)-2K_1(x)^{2}}{2K_0(x)^2},
\]
and we can use the inequality
\(K_0(x)^2+K_2(x)K_0(x)-2K_1(x)^{2}>0\) (which can be verified
numerically).
\newline\indent
Thus, if \(l<1\), we have
\(\frac{V'(\frac{r}{l})}{V(\frac{r}{l})}\geq \frac{V'(r)}{V(r)}\).
Since \(V(r)<0\), \(V'(r)>0\), this implies
\({V'(\frac{r}{l})}{V(r)}-{V'(r)}{V(\frac{r}{l})}\geq 0\), and
therefore
\[
\begin{split}
h'(r) & = \frac{C}{l(V'(r))^{2}}\big(
\tfrac{1}{l}V'(r)V(\tfrac{r}{l})-V'(\tfrac{r}{l})V(r)\big)
\\
& <\frac{C}{l(V'(r))^{2}}\big(
V'(r)V(\tfrac{r}{l})-V'(\tfrac{r}{l})V(r)\big)\leq 0.
\end{split}
\]
Similarly, if \(l>1\), we obtain \(h'(r)>0\).
\newline\indent
Further, it is directly checked that
$$
\lim\limits_{r\to0+} h(r) = Cl^{n-2}\quad \mbox{and}\quad
\lim\limits_{r\to\infty} h(r) = \left\{\begin{array}{lc} 0 &
\mbox{ if } l<1\\
+\infty & \mbox{ if } l>1
\end{array}\right. .
$$
Thus, the equation \(h(r)=1\) has no solution in the cases
\(Cl^{n-2}<1,\;l<1\) or \(Cl^{n-2}>1, l>1\) and a unique positive
solution in the cases \(Cl^{n-2}<1,\;l>1\) or \(Cl^{n-2}>1, l<1\).
Recalling that \(U'(r)=V'(r)(1-h(r))\), we see that of the last
two cases, the former corresponds to a local maximum of \(U(r)\)
and the latter to a local minimum. Since the minimum is unique,
then it is a global minimum.
\newline\indent
Concerning the second assertion, by construction we have
$$
\int_{\rrn}V(|x|)\dd x=-1 \quad \mbox{for all } n\,,
$$
since $V$ is the fundamental solution of $\Delta u - u =
\delta_{0}$ as stated in Definition \ref{defv}. Therefore, we get
\begin{equation*}
\int_{\rrn}U(|x|)\dd x = \int_{\rrn}V(|x|)-CV(|x|/l) \dd x
=-1+Cl^{n},
\end{equation*}
which is negative for $Cl^{n}<1$, and thus $U$ is catastrophic
(see \cite{rue}, p.\@ 37).

\section{Acknowledgements}
JAC acknowledges partial support by MICINN project, reference
MICINN MTM2011-27739-C04-02, by GRC 2009 SGR 345 by the
Generalitat de Catalunya, and by the Engineering and Physical
Sciences Research Council grant number EP/K008404/1. JAC also
acknowledges support from the Royal Society through a Wolfson
Research Merit Award. We also acknowledge the Isaac Newton
Institute for the Mathematical Sciences, where part of this work
was accomplished. We are grateful to Yanghong Huang for several
useful comments.


\end{document}